\documentclass[12pt]{article}
\usepackage[margin=1in]{geometry}
                    
\usepackage{cite}
\usepackage{amsmath,amssymb,amsfonts}
\usepackage{algorithm}
\usepackage{algpseudocode}
\usepackage{caption}
\usepackage{subcaption}
\usepackage{graphicx}
\usepackage{textcomp}
\usepackage{xcolor}
\usepackage{amsmath}
\usepackage[hidelinks]{hyperref}
\def\BibTeX{{\rm B\kern-.05em{\sc i\kern-.025em b}\kern-.08em
    T\kern-.1667em\lower.7ex\hbox{E}\kern-.125emX}}

\title{\LARGE \bf
A Dynamic Unmanned Aerial Vehicle Routing Framework for Urban Traffic Monitoring\\}

\author{Yumeng Bai and Yiheng Feng
\thanks{This research is supported in part by the U.S. National Science Foundation (NSF) through Grant CPS 2038215, and U.S. Department of Transportation (USDOT) Region 5 University Transportation Center (Center for Connected and Automated Transportation). \textit{(Corresponding author: Yiheng Feng.)}}
\thanks{The authors are with the Lyles School of Civil and Construction Engineering, Purdue University, West Lafayette, IN 47907 USA (e-mail: bai145@purdue.edu; feng333@purdue.edu).}
}

\begin{document}

\maketitle
\thispagestyle{empty}
\pagestyle{empty}

\begin{abstract}
Unmanned Aerial Vehicles (UAVs) have great potential in urban traffic monitoring due to their rapid speed, cost-effectiveness, and extensive field-of-view, while being unconstrained by traffic congestion. 
However, their limited flight duration presents critical challenges in sustainable recharging strategies and efficient route planning in long-term monitoring tasks.
Additionally, existing approaches for long-term monitoring often neglect the dynamic and evolving nature of urban traffic networks.
In this study, we introduce a novel dynamic UAV routing framework for long-term, network-wide urban traffic monitoring, leveraging existing ground vehicles (GVs), such as public transit buses, as mobile charging stations without disrupting their operations. 
To address the complexity of long-term monitoring scenarios involving multiple flights, we decompose the problem into manageable single-flight tasks, in which each flight is modeled as a Team Arc Orienteering Problem with Decreasing Profits (TAOP-DP) with the objective to collectively maximize the spatiotemporal network coverage. Between flights, we adaptively update the edge weights to incorporate real-time traffic changes and revisit intervals.
We validate our framework through extensive microscopic simulations in a modified Sioux Falls network under various scenarios, including static and changing demand levels and varying historical information availability.
Comparative results demonstrate that our model outperforms three baseline approaches, especially when historical information is incomplete or absent.
Moreover, we show that our monitoring framework can capture network-wide traffic trends and construct accurate Macroscopic Fundamental Diagrams (MFDs).
These findings demonstrate the effectiveness of the proposed dynamic UAV routing framework, underscoring its suitability for efficient and reliable long-term traffic monitoring. Our approach's adaptability and high accuracy in capturing the MFD highlight its potential in network-wide traffic control and management applications.
\end{abstract}

\newpage

\section{Introduction}
\label{sec:intro}

Unmanned Aerial Vehicles (UAVs) have gained increasing interest in recent years due to their versatility, rapid speeds, cost-effectiveness, and extensive field-of-view. These characteristics make UAVs vital for future smart cities, broadening their roles in areas such as package delivery~\cite{uav-delivery}, road inspection~\cite{multi-depot-recharge}, crowd surveillance~\cite{surveillance} and traffic monitoring~\cite{pNEUMA, space-time-uav-monitor, uav-monitor-1, uav-monitor-2, uav-monitor-3, uav-monitor-4, urban-monitoring}.
As urban areas grow and traffic becomes more complex, there is an urgent need for comprehensive and long-term traffic monitoring systems that can capture network-level traffic conditions and congestion propagation and dissipation to ensure efficient and safe transportation network operations~\cite{pNEUMA}.

This increasing need has driven recent studies into utilizing vehicles, including UAVs, connected vehicles (CVs) and connected and automated vehicles (CAVs) as mobile sensors for traffic monitoring~\cite{space-time-uav-monitor, cv-monitor-1, cv-monitor-2}. 
Compared to ground vehicles, UAVs offer unique advantages. CV-based traffic monitoring usually requires a critical penetration rate that may not be achieved in the near future. Although CAV-based (i.e., cooperative perception) traffic monitoring can significantly reduce the required penetration rate, the data collection efficiency and quality is greatly impacted by traffic conditions and sensor's detection capabilities and locations on the vehicle. 
In contrast, UAVs are unconstrained by road networks or congestion, and can reach locations inaccessible to fixed sensors. Moreover, data captured by UAVs provide a bird-eye-view of traffic flow with extended detection range and limited occlusions.
This capability enables UAVs to gather dynamic, timely data, significantly enhancing traffic monitoring in urban environments~\cite{surveillance, pNEUMA}.

Given these advantages, extensive research has explored UAV-based traffic monitoring from various aspects~\cite{pNEUMA, space-time-uav-monitor, uav-monitor-1,uav-monitor-2,uav-monitor-3}.
For instance, a field experiment in \cite{pNEUMA} employed a swarm of UAVs hovering over urban road networks, validating their effectiveness in collecting comprehensive traffic data in a multi-modal congested environment.
Additionally, \cite{uav-monitor-3} developed a framework to simultaneously detect and count vehicles from UAV-captured images.
The approach proposed in \cite{uav-monitor-2} optimizes computation offloading decisions for UAVs, enabling efficient data processing for real-time operations.
These studies collectively lay the foundation for UAV-based data collection and processing. 
However, most existing studies focus on short-term operations at static locations and overlook the needs of long-term tasks that require multiple flights covering different areas.

Achieving long-term and dynamic traffic monitoring with UAVs presents significant challenges due to their limited flight durations caused by energy constraints. Consequently, developing effective recharging strategies and efficiently planning UAVs' routes are critical~\cite{surveillance, pNEUMA, multi-depot-recharge, urban-monitoring}.
Two main strategies have been proposed to address UAV recharging requirements. 
The first involves using multiple fixed depots where UAVs can recharge or swap batteries~\cite{multi-depot-recharge, urban-monitoring}. While this approach reduces computational costs, the static locations of depots limit the monitoring range of UAVs to specific areas, requiring additional depots to cover larger regions. 
Alternatively, the cooperation of ground and aerial vehicles using ground vehicles (GVs) as mobile charging stations has been increasingly studied for various tasks, such as using trucks and drones cooperatively for delivery~\cite{drone-truck} or employing drones and ground vehicles collectively for surveillance~\cite{uav-ugv-recharge, air-ground}.
This method allows UAVs to recharge while in transit, but dedicating GVs specifically for UAV recharging significantly increases the overall costs and logistical complexity. 
Previous studies have also explored UAV recharging solutions utilizing existing GVs in the network with fixed routes \cite{fixed-moving-gv,drone-bus-delivery}. However, \cite{fixed-moving-gv} focused on minimizing mission completion time, and the framework in \cite{drone-bus-delivery} targeted delivery tasks, neither of which are designed for continuous traffic monitoring task.

Methodologically, most of the existing literature on drone or air-ground joint routing formulates the problem as a Vehicle Routing Problem (VRP), aiming to minimize the total travel time or distance to visit all targets \cite{drone-truck, air-ground, drone-routing-1, drone-routing-2}. 
In contrast, the urban traffic monitoring problem targets maximizing spatiotemporal coverage using multiple UAVs over a given time period (e.g., an entire day with demand fluctuations), which is more close to a Team Orienteering Problem (TOP) \cite{TOP}. 
Moreover, current research on persistent monitoring often involves repetitive cycles of visiting all targets, minimizing time between revisits, or imposing a revisit period constraint~\cite{persistent-1, persistent-2, persistent-3}. 
Such static approach does not incorporate the dynamic nature of urban traffic networks, potentially leading to inefficiencies.
Other studies on TOP with multiple visits or decaying profits typically employ an arbitrarily assumed fixed decay factor~\cite{TOP-DP, TOP-DP-2, TOP-DP-3, TOP-multi-visit, TOP-multi-visit-2, OTOP-RV}, disregarding the dynamic evolution of traffic conditions in both the spatial and temporal domains.

\begin{figure*}[!t]
  \begin{center}
    \includegraphics[width=\textwidth]{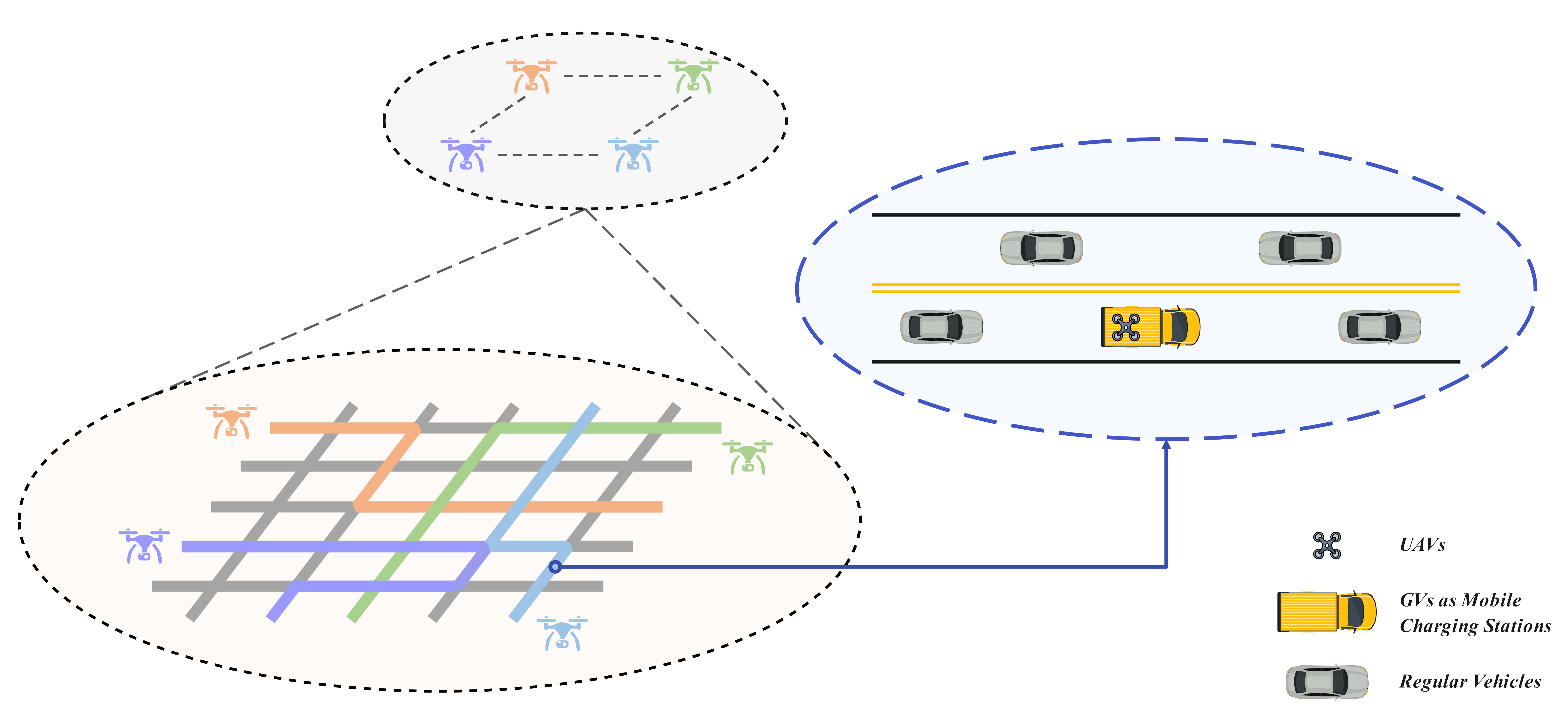}
  \end{center}
  \caption{Conceptual overview of the proposed UAV-Based traffic monitoring framework.}
  \label{fig:concept}
\end{figure*}

To overcome these limitations, we propose a novel dynamic UAV routing framework for long-term traffic monitoring that leverages existing fixed-route GVs, such as public transit buses, as dynamic charging stations, to enhance cost-efficiency and expand the monitoring range, as illustrated in Figure~\ref{fig:concept}.
Unlike previous cooperative solutions requiring GVs to stop or detour to meet with UAVs, our framework uses predefined GV routes as constraints for the UAV routing problem. This approach maximizes the utilization of the current fleet without impacting their operations.
We formulate the routing problem for each single flight with multiple UAVs as an TOP problem, which seeks to collectively maximize the total score obtained by visiting a subset of edges.
To adapt to the fluctuations of the traffic demand and improve monitoring efficiency, we further incorporate a design with decreasing profits over time to prioritize high-priority edges, resulting in an arc-based Team Arc Orienteering Problem with Decreasing Profits (TAOP-DP)~\cite{TOP-DP, TOP-DP-2, TOP-DP-3}. Then a novel heuristic algorithm which decomposes the long-term monitoring task into multiple single flights is developed to reduce the computation complexity and incorporate dynamic traffic information collected in previous flights and arc revisit intervals.

Furthermore, existing UAV-based traffic monitoring studies primarily focus on collecting link-level data such as vehicle detection, tracking, counts, and speeds~\cite{pNEUMA, uav-monitor-3, uav-monitor-4}.
While the effectiveness and accuracy of UAV-based traffic monitoring at the link level has been extensively explored and verified, few studies have looked into network-level properties and performance. To fill this research gap, we show that our dynamic monitoring framework can capture network-wide traffic trends and construct an accurate Macroscopic Fundamental Diagram (MFD). To the best of our knowledge, this is among the earliest works that utilize UAV data to generate MFD.

The main contributions of this study can be summarized as follows:
\begin{enumerate}
    \item We present a dynamic UAV routing framework for long-term traffic monitoring that leverages existing GVs as mobile charging stations for UAVs, significantly enhancing cost-efficiency and expanding the monitoring coverage without disrupting GV operations.
    \item We formulate the single-flight UAV routing problem as a TAOP-DP to optimize spatiotemporal network coverage and develop a novel heuristic algorithm that adaptively updates edge priorities between flights based on real-time traffic data and information decay.
    \item We conduct a microscopic simulation to evaluate the performance of UAV routing algorithms, demonstrating the robustness and applicability of our approach under various scenarios and establishing a basis for future vehicle-level cooperation between GVs and UAVs.
    \item  We expand the scope of UAV-based traffic monitoring from link-level accuracy to network-level performance. The promising results validate the accuracy and effectiveness of our method in MFD generation and provide a foundation for real-time traffic management strategies informed by MFD analysis.
\end{enumerate}

The remainder of the paper is organized as follows. Section~\ref{sec:framework} introduces the overall framework. Section~\ref{sec:methodology} describes the methodologies of the proposed dynamic UAV routing framework. Section~\ref{sec:case} presents a case study on the modified Sioux Falls network and analyzes the performance of the framework in various scenarios. Section~\ref{sec:mfd} evaluates the framework's ability to capture network-level traffic patterns through MFD analysis. Section~\ref{sec:discussion} discusses the applications and extensions of the proposed framework. Finally, Section~\ref{sec:conclusion} provides conclusions.

\section{FRAMEWORK OVERVIEW}
\label{sec:framework}
\begin{figure*}[!t]
  \begin{center}
    \includegraphics[width=\textwidth]{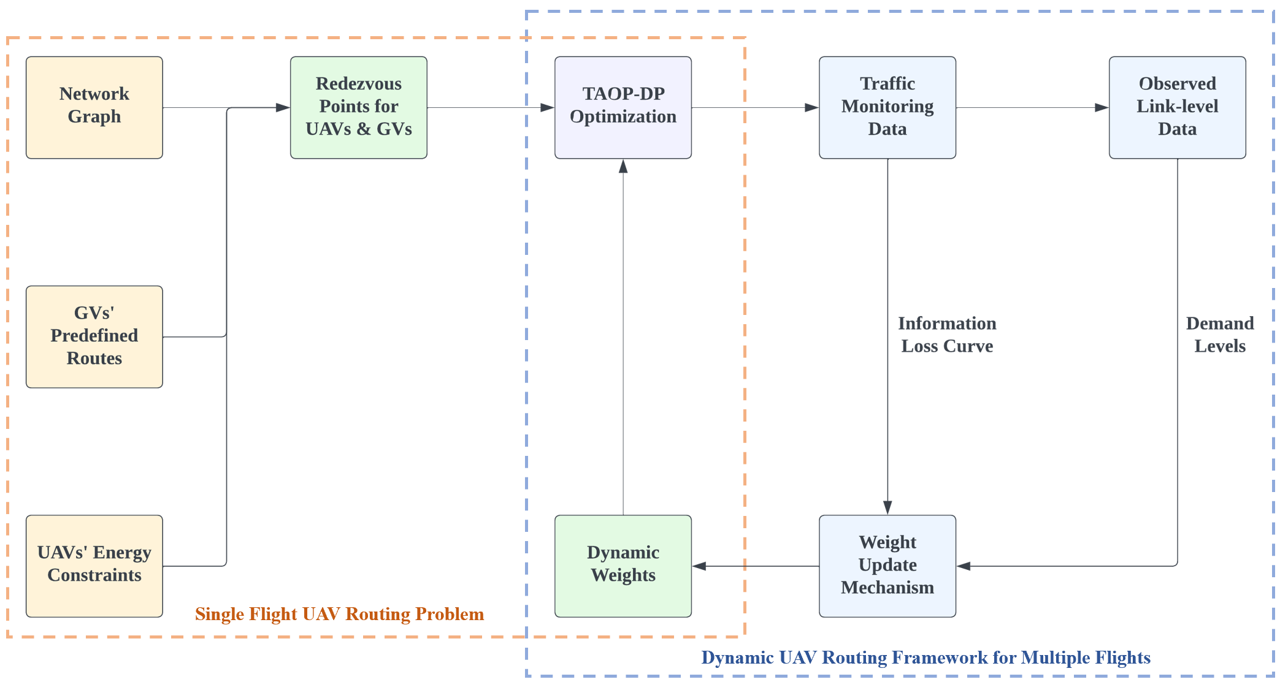}
  \end{center}
  \caption{Overview of the dynamic UAV routing framework.}
  \label{fig:overall}
\end{figure*}

We propose an efficient and adaptive solution for long-term urban traffic monitoring that effectively captures both link-level and network-level traffic patterns, including the MFD. Figure \ref{fig:overall} illustrates the overview of our framework. We model each multi-UAV single-flight routing problem as a TAOP-DP problem, with UAV's energy and rendezvous points from GVs as constraints. These rendezvous points are GV locations at the end of the flight duration, serving as mobile charging stations. We assume that the GVs' routes are known (e.g., transit buses) so that the rendezvous points can be estimated before each flight. In the experiments, we demonstrate that our framework can be easily adapted to support adaptive rerouting, where the estimation of the rendezvous is not accurate due to uncertainties caused by congestion, dwell time variation etc. We further assume that charging is performed by swapping the battery so that charging time is not considered. This assumption is supported by developments such as the automated battery swapping system presented in~\cite{battery-swap}.

For long-term traffic monitoring involving multiple flights, we develop a heuristic algorithm that decomposes the task into several single-flight problems, integrating real-time traffic data collected from previous flights and revisit intervals. 
Inspired by \cite{info-TSE}, we introduce an information loss curve calibrated from traffic data to dynamically capture how information loss changes over revisit intervals and use the curve to model the revisit value decay. 
After executing optimal routing plans and data collection for each single flight, we implement a weight update mechanism for the next flight to incorporate real-time monitoring results. This mechanism combines changes in traffic demand based on the latest observed data and information decay based on revisit intervals using the information loss curve.

To validate our framework, we conduct a microscopic simulation study on the modified Sioux Falls network under both static and changing demand scenarios. The results demonstrate the robustness and effectiveness of our dynamic UAV routing framework under various conditions compared to several baselines, especially when prior information (e.g., historical link flow data) is incomplete or unavailable. Additionally, our framework accurately captures the MFD at the network-level, highlighting its potential for future dynamic network traffic management.

\section{METHODOLOGY} 
\label{sec:methodology}
In this section, we present the methodology of the dynamic UAV routing framework for urban traffic monitoring. We begin by modeling the multi-UAV single-flight routing problem as an arc-based Team Arc Orienteering Problem with Decreasing Profits (TAOP-DP).
Next, we propose a heuristic algorithm that solves the dynamic UAV routing problem with multiple flights by adaptively combining multiple TAOP-DP problems, and continuously updating edge weights base on routing decisions and observed data to capture real-time traffic changes. Finally, we introduce the information loss curve used in the weight update mechanism and describe the entire algorithm in detail.
Main notations are summarized in Table~\ref{tab:parameters}.

\begin{table}[!t]
\centering
\caption{Variables and Notations}
\label{tab:parameters}
\resizebox{\columnwidth}{!}{%
\begin{tabular}{ll}
\hline
\textbf{Variables} & \textbf{Descriptions} \\ \hline
$K$ & Set of UAVs.\\
$V$ & Set of vertices in the directed graph $G(V,E)$.\\
$E$ & Set of edges in the directed graph $G(V,E)$.\\
$H$ & Set of flights.\\
$M$ & A large positive constant. \\ 
$c$ & Decay constant. \\ 
$D_{i,j}^h$ & Assigned demand level of edge $(i,j)$ in flight $h$. \\
$w_{i,j}^h$ & Assigned weight of edge $(i,j)$ in flight $h$. \\
$d_{i,j}^h$ & Decay ratio for the weight of edge $(i,j)$ in flight $h$, equals to $w_{i,j}^h/c$. \\
$t_{i,j}$ & Travel time from vertex $i$ to vertex $j$. \\
$T_h$ & Time budget for each agent in flight $h$. \\ 
$s_k$ & Start vertex for agent $k$. \\ 
$e_k$ & End vertex for agent $k$. \\ 
$x_{i,j,k}$ & Binary variable. 1 if agent $k$ travels from vertex $i$ to vertex $j$, and 0 otherwise.\\
$a_{i,k}$ & Non-negative continuous variable. The visiting time of vertex $i$ by agent $k$.\\ 
$z_{i,j,k}$ & Non-negative continuous variable. The departure time of edge $(i,j)$ by agent $k$.\\ 
$I(t)$ & Information loss function that characterizes the information loss change over time.\\ \hline
\end{tabular}%
}
\end{table}

\subsection{Single-flight UAV Routing Problem Formulation}

To begin, we formulate the UAV routing problem for each single flight as a Mixed Integer Linear Programming (MILP) problem. To leverage moving GVs as recharging stations, several key factors are integrated into our model: the network graph, the GVs' predefined routes, and the UAVs' energy constraints. These factors collectively determine the spatial and temporal rendezvous points where UAVs and GVs can meet for recharging (e.g., battery exchange). These rendezvous points serve as constraints in the UAVs' routing problem, and may change due to changing traffic conditions, requiring re-generation of the UAVs routes. 
Additionally, we incorporate varying edge weights based on traffic demand, prioritizing high-demand areas in UAV routing. 
The overall structure of the single-flight UAV routing problem is shown in the left part of Figure \ref{fig:overall}.

Within each flight, UAVs aim to jointly maximize spatiotemporal coverage of the traffic network, considering varying edge weights. The objective function maximizes the total score obtained by visiting a subset of edges with teams of UAVs within a given time duration. 
Due to the characteristics of traffic monitoring, we prioritize visiting high-priority edges as soon as possible. High-priority edges are determined by their 1) traffic demand levels and 2) elapsed time since the last visit. To achieve this, we introduce decreasing profits, formulating the problem as a TAOP-DP, which is an arc-based variant of \cite{TOP-DP}.
A distinctive feature of our model is the decay ratio \( d_{i,j}^h \), which is set in direct proportion to the dynamic weight \( w_{i,j}^h \). Unlike the fixed decay ratio in existing literature, the dynamic weight ensures higher-priority edges are visited first. The mathematical formulation (\textbf{P1}) is presented as follows:

\begin{align}
\max \quad & \sum_{(i,j)\in E}\sum_{k\in K} ((w_{i,j}^h + w_{j,i}^h) x_{i,j,k} - (d_{i,j}^h + d_{j,i}^h) z_{i,j,k}) \\
\text{s.t.} \quad & \nonumber\\
& \sum_{i \in V} x_{s_k,i,k} = 1 &\quad \forall k \in K \\
& \sum_{i \in V} x_{i,e_k,k} = 1 &\quad \forall k \in K\\
& \sum_{(i,l) \in E} x_{i,l,k} = \sum_{(l,j) \in E} x_{l,j,k} &\quad \forall l \in V \backslash \{s_k, e_k\}, \forall k \in K\\
& \sum_{k \in K} (x_{i,j,k} + x_{j,i,k})  \leq 1 & \quad \forall (i, j) \in E\\
& a_{i,k} \geq 0 &\quad \forall i \in V, \forall k \in K\\
& a_{i,k} + t_{i,j} - a_{j,k} \leq M (1 - x_{i,j,k}) &\quad \forall (i, j) \in E, \forall k \in K\\
& \sum_{(i,j) \in E} x_{i,j,k} t_{i,j}\leq T_h &\quad  \forall k \in K \\
& x_{i,j,k} \in \{0,1\} &\quad  \forall  (i, j) \in E, \forall k \in K \\
& z_{i,j,k} \geq 0 &\quad \forall (i, j) \in E, \forall k \in K\\
& z_{i,j,k} \geq a_{i,k} + t_{i,j} - M (1 - x_{i,j,k}) &\quad   \forall  (i, j) \in E, \forall k \in K \\
& z_{i,j,k} \leq  M x_{i,j,k} &\quad   \forall  (i, j) \in E, \forall k \in K
\end{align}\\

Due to the nature of UAV-based traffic monitoring, we assume that when an edge $(i,j)$ is visited, its reverse edge $(j,i)$ is also visited.
The objective function (1) maximizes the total weighted scores within the current flight $h$ collected by all UAVs. The first term (primary objective) selects the edges to be visited during the flight, while the second term (secondary objective) optimizes the visiting sequence of those selected edges. In Equation (1), $d^h_{i,j}=w^h_{i,j}/c$, where $c$ is a decay constant and is selected to be half of the flight duration (i.e., 15) in this study to ensure that the primary objective takes precedence over the secondary objective. Constraints (2)-(3) guarantee that each UAV departs from its start point $s_k$ and returns to its end point $e_k$ (i.e., the rendezvous point with the GV). Constraint (4) ensures the flow conservation for each vertex. Constraint (5) guarantees that each edge, and its reverse edge, is visited at most once within the current flight $h$. Constraints (6), (9) and (10) define the nature of variables $a_{i,k}$ (non negative), $x_{i,j,k}$ (binary) and $z_{i,j,k}$ (non negative). Constraint (7) prevents the generation of subtours. Constraint (8) ensures that all UAVs return to their end points within the time budget $T_h$. Constraint (11)-(12) ensure that the departure time at each edge for UAV $k$, $z_{i,j,k}$ accounts for the visiting time of the preceding vertex plus the travel time if the edge is traversed, otherwise it is 0.

\subsection{Dynamic UAV Routing Framework for Multiple Flights}
When considering traffic monitoring over an extended period involving multiple flights, the complexity of the problem increases significantly. Unlike the single-flight routing problem, where revisit is not considered due to UAV energy constraints, multiple flights with longer operation horizons allow revisiting edges when necessary. It is crucial for UAVs to explore previously unvisited edges or those visited long time ago to capture evolving traffic conditions. To address these challenges, we propose a dynamic routing framework for multiple flights, as illustrated in the right part of Figure \ref{fig:overall}. 

Our approach decomposes the overall problem into multiple single-flight problems. After executing the optimal routing plan for each flight and collecting traffic data, we update the edge weights for the next flight based on two factors: the time since each edge was last visited (revisit interval) and the current demand level derived from the observed data. These updated weights guide UAVs to prioritize edges with higher traffic demand or those not visited for a long time. This decomposition offers two main benefits: 1) it greatly reduces the computational cost; and 2) the weight update enables the routing framework to respond to real-time traffic changes, enhancing monitoring efficiency.

The weight update mechanism consists of two key components. First, the importance of edges may change due to evolving traffic conditions, such as accidents or congestion in specific locations. We account for this by classifying observed link-level flow data into different demand levels (i.e., each demand level is assigned one weight).
Second, the mechanism considers the impact of revisit intervals (i.e., time elapsed since the edge was last visited). Intuitively, shorter revisit intervals generally yield less valuable new information due to limited changes in traffic conditions. By incorporating both real-time traffic conditions and the information loss of frequent visits, our routing framework dynamically maximizes monitoring efforts.

\subsubsection{Information Loss Curve}
We propose an information loss curve $I(t)$ to model how the value of monitoring information changes over revisit intervals. The process to capture, fit, and calibrate $I(t)$ is shown in Algorithm \ref{alg:info}. Given a road network with $|E|$ edges, information loss is determined by the errors in the latest observed link-level data. 
To calculate the observation error, we randomly select a specified number of edges and classify them as observed, while the remaining edges are set as unobserved.
For the observed edges, their traffic states are directly collected. For unobserved edges, we estimate their current states using the most recently observed data, or assume zero if they have never been observed. The estimation error is then computed by comparing these estimates to the ground truth data.
We repeat this process ten times for each specified number of edges and take the average results to accommodate the randomness in the edge selection process. By varying the number of selected edges from one to the total number of edges, we calculate the corresponding revisit intervals. Finally, by plotting the estimation errors against the revisit intervals, we generate and fit the information loss curve $I(t)$, which quantifies how information loss evolves over time. A concrete example is provided in the Case Study section.

\begin{algorithm}
\caption{Pseudocode for Information Loss Curve Calibration.} 
\label{alg:info}
\begin{algorithmic}[1]
\State \textbf{Input:} Ground truth data $G$ and time interval duration $\Delta t$
\State \textbf{Output:} Information loss curve $I(t)$

\For{$S = 1$ to $|E|$}
    \State Initialize $Sum\_Error \gets 0$
    \For{$i = $1 to 10}
        \For {each time interval in $G$}
            \State Randomly select $S$ edges to observe.
            \State Record traffic states of observed edges.
            \State Estimate traffic states for unobserved edges using the latest observed data or 0.
        \EndFor
        \State $Current\_Error \gets$ Compute estimation error with ground truth data
        \State $Sum\_Error \gets Sum\_Error + Current\_Error$
    \EndFor
    \State $Avg\_Error[S] \gets \frac{Sum\_Error}{10}$
    \State $Revisit\_Interval[S] \gets \frac{|E|}{S} \cdot \Delta t$
\EndFor
\State \textbf{Plot:} $Avg\_Error$ vs. $Revisit\_Interval$
\State \textbf{Fit Curve:} Derive $I(t)$ to represent information loss over revisit intervals.
\end{algorithmic}
\end{algorithm}

\subsubsection{Multi-Flight Heuristic Algorithm.}
Building upon the information loss curve $I(t)$, we propose a heuristic algorithm to solve the multi-flight dynamic UAV routing problem, as shown in Algorithm \ref{alg:decomp_algorithm}.

\begin{algorithm}
\caption{Pseudocode for Multi-Flight Heuristic Algorithm} \label{alg:decomp_algorithm}
\begin{algorithmic}[1]

\State \textbf{Step 0: Initialization}
\State Initialize the weight and demand level of each edge $(i,j)$ as $D_{i,j,0}$.
\State Initialize the revisit interval $v_{i,j}$ of each edge $(i,j)$ as 0.
\State Define the information loss function $I(t)$.

\For{each flight $h \in H$}
    \State \textbf{Step 1: Solve \textbf{P1} for the current flight $h$.}
    \If {rendezvous points change or flight $h$ completes}
        \State \textbf{Step 2: Update edge weights based on executed routes and observed data.}
        \For{each edge $(i,j)$}
            \If{edge $(i,j)$ is visited in flight $h$}
                \State  $v_{i,j} \gets T_h - \sum_{k \in K}(z_{i,j,k}+z_{j,i,k})$
            \Else
                \State $v_{i,j} \gets v_{i,j} + T_h$
            \EndIf
            \State $D_{i,j,h} \gets\text{Observed link-level data}$
            \State $w_{i,j,h+1} \gets D_{i,j,h} \cdot I(v_{i,j})$
            \State $d_{i,j,h+1} \gets \frac{w_{i,j,h+1}}{c}$
        \EndFor
        \State \textbf{Step 3: Update constraints for the next flight.}
        \For{each agent $k$}
            \State $StartPoint_{k,h+1} \gets CurrentPosition_{k}$
            \State $EndPoint_{k,h+1} \gets \text{Next rendezvous point}$
            \State $T_{h+1} \gets$ remaining flight duration if not zero, otherwise new flight duration.
        \EndFor
    \EndIf
\EndFor

\end{algorithmic}
\end{algorithm}

We first initialize the demand levels and revisit intervals for all edges based on the available information. We then solve the TAOP-DP problem for the current flight $h$ as described in the single-flight scenario. 
During the execution of the routing plan, if the rendezvous points change or the current flight completes, we proceed to update the edges weights based on the executed routes and newly observed data. The weight update process is as follows: first, we update the revisit interval $v_{i,j}$ for each edge based on whether it was visited during the flight. Next, we update demand levels $D_{i,j,h}$ using the observed data. The new weight for the next flight is calculated by multiplying the demand level with the information loss $I(v_{i,j})$, determined by the information loss curve. The corresponding decay ratio $d_{i,j,h+1}$ is computed as $\frac{w_{i,j,h+1}}{c}$.
Additionally, we update the constraints for the next flight by setting the current position of each UAV as the new starting point and assigning the next rendezvous point as the endpoint. The duration of the next flight $T_{h+1}$ is updated to either the remaining flight duration, if not zero, or a new flight duration.

By integrating the information loss curve and adapting edge weights in real-time, our heuristic algorithm effectively addresses the computational challenges of multi-flight UAV routing while ensuring responsiveness to dynamic traffic conditions. In the following sections, we validate the effectiveness of our adaptive framework through microscopic simulations under both static and varying demand level scenarios, assessing its performance at the link-level and its capability to capture network-level traffic patterns using MFD.

\section{CASE STUDY}
\label{sec:case}

In this section, we conduct a simulation study on the modified Sioux Falls network to evaluate the effectiveness of our proposed dynamic UAV routing framework. 
We begin by describing the simulation environment setup and generating the corresponding information loss curve $I(t)$. Following this, we detail the experimental setup and introduce three baseline models for comparison.
The performance of our model is first analyzed under static demand conditions and benchmarked against the baseline models.
Then it is extended with a sensitivity analysis to evaluate the impact of  varying levels of historical information availability. 
Subsequently, we examine a changing demand scenario that reflects real-world daily traffic dynamics.

\subsection{Simulation Environment Setup}
To demonstrate our model, we employ the widely used Sioux Falls network, a standard benchmark in traffic network modeling studies~\cite{sioux-eq}. The network is extracted from the Sioux Falls city in the USA, containing 24 vertices and 76 edges, with edge lengths ranging from 0.4 to 3.4 miles. We choose this network because it is a representative example of a small urban network, making it ideal for our simulation study. To generate ground truth data for constructing the information loss curve and evaluating our model, we build this network in SUMO microscopic traffic simulator~\cite{sumo}.

Travel demand dataset from~\cite{sioux-data} is used for vehicle generation after appropriate scaling. Vehicle trips are generated using the Origin-Destination (OD) demand from this dataset and follow a uniform distribution. These trips are then assigned to routes based on dynamic user equilibrium, implemented using SUMO's built-in functions~\cite{sumo-dua}. Each vertex is considered as a signalized intersection and fixed-time traffic signal is applied at each intersection. Split phasing is applied, with each phase duration proportional to the average demand, and the total cycle length is set to 100 seconds for all intersections.

\begin{figure}[!ht]
  \begin{center}
    \includegraphics[width=0.8\textwidth]{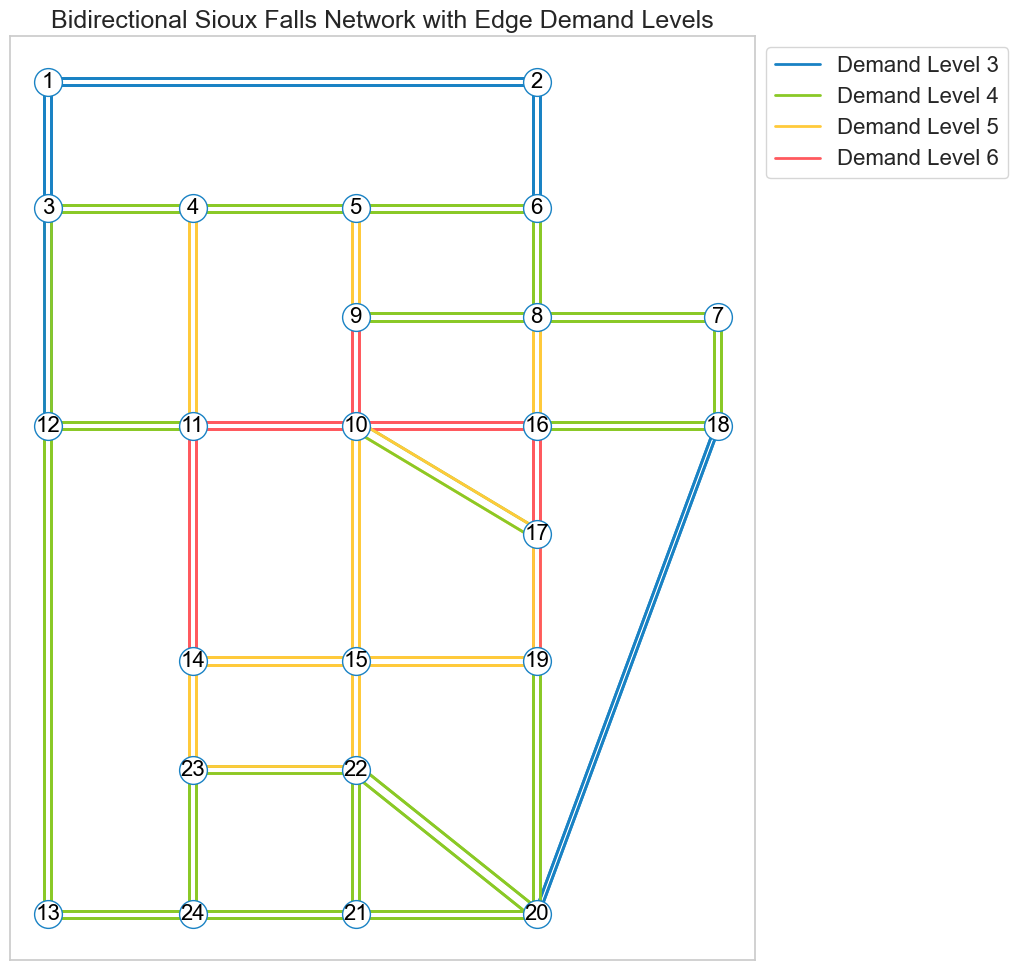}
  \end{center}
  \caption{Bidirectional Sioux Falls network with edge demand levels.}
  \label{fig:network}
\end{figure}

We run the simulation for 102 hours and use the latter 100 hours for data collection. During this period, we record traffic flow data at ten-minute intervals, defined as the number of vehicles traversing each edge within each interval. The data served as the ground truth for the information loss curve generation. The average flow value is aggregated per ten-minute interval and analyzed over the entire 100-hour simulation to assign varying demand levels to different edges, as shown in Equation \ref{eq:demand}. These demand levels represent the relative weight of each edge (i.e., $D$) in the optimization model. The demand level distribution of the Sioux Falls network is shown in Figure \ref{fig:network}. Note that in real world implementation, the average flow rate (or an approximation) can be obtained from other traffic data such as Annual Average Daily Traffic (AADT).

\begin{equation}
    D=
    \begin{cases} 
    3 & \text{if } \text{average flow} \leq 15 \\
    4 & \text{if } 15 < \text{average flow} \leq 30 \\
    5 & \text{if } 30 < \text{average flow} \leq 45 \\
    6 & \text{if } \text{average flow} > 45 
    \end{cases}
\label{eq:demand}
\end{equation}

\subsection{Information Loss Curve Generation}
\begin{figure}[!ht]
  \begin{center}
    \includegraphics[width=0.8\textwidth]{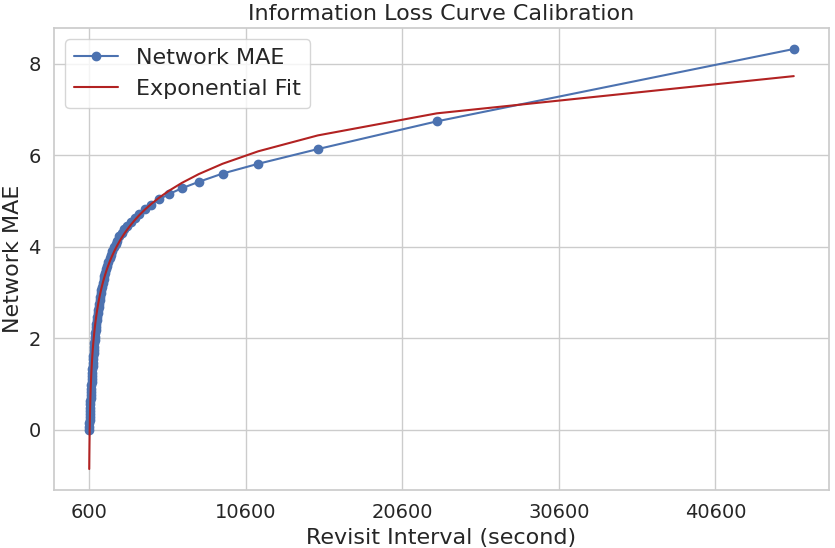}
  \end{center}
  \caption{Information loss curve calibration result.}
  \label{fig:info}
\end{figure}

Using the ground truth data collected above, we select Mean Absolute Error (MAE) of flow as our information error metric. The information loss curve is generated by applying Algorithm \ref{alg:info}.
The resulting estimation error over revisit intervals, along with the fitted curve, is presented in Figure \ref{fig:info}.
After fitting the curve using an exponential function, we normalize it and obtain the following equation:
\begin{equation}
    I(t) = (1.1531 \cdot \ln(x - 573.8506) - 4.6229)/2
    \label{eq:info}
\end{equation}
The $R^2$ is 0.98, indicating an excellent fit.
The shape of the curve matches the typical trend of information loss (i.e., exponential delay over time) where the error increases rapidly with longer revisit intervals at the beginning and gradually stabilizes over extended periods.

\subsection{Experimental Design}
For each single flight, the flight duration is set to 30 minutes, and the UAV speed is 1000 meters per minute, assumed to be constant. The travel time of each edge is calculated by dividing the edge length by the UAV speed.

Due to the network size, flight duration, and UAV speed, we consider two UAVs in the experiment. The four rendezvous points (vertices) for UAV 1 is 4-22-6-24. After the UAV reaches vertex 24, it will come back to vertex 4 and repeat the loop. Similarly, the rendezvous points for UAV 2 are 22-4-24-6 and repeat afterwards. This setting mimics GV operations for mobile charging, where these rendezvous points can be considered as depots or start/end locations of GV routes (e.g., bus stops or depots). 
For each multiple-flight experiment, we perform 20 iterations, equivalent to 10 hours of traffic monitoring.

\subsection{Baseline Models for Comparison}
In this section, we introduce three baseline models for comparison with our dynamic framework.
\subsubsection{Baseline with Linear Information Loss Curve (\textbf{Linear Info})}
We consider a baseline model with a linear information curve where, instead of using our fitted information loss curve, a linear approximation is applied. Given that a 10-minute time step is used when generating the information loss curve, we assume that traffic conditions remain unchanged within each time step (i.e., $I(t) = 0$ for $t \leq 600$ seconds). We further assume that the information loss increases linearly, with a slope proportional to the demand level divided by the flight duration. The equation for the linear curve is defined as follows:
\begin{equation}
I(t) =
\begin{cases} 
0 & \text{if } t \leq 600 \ \text{seconds} \\
\frac{D(t - 600)}{T_h - 600} & \text{if } t > 600 \ \text{seconds}
\end{cases}
\end{equation}

\begin{figure}[!ht]
  \begin{center}
    \includegraphics[width=\textwidth]{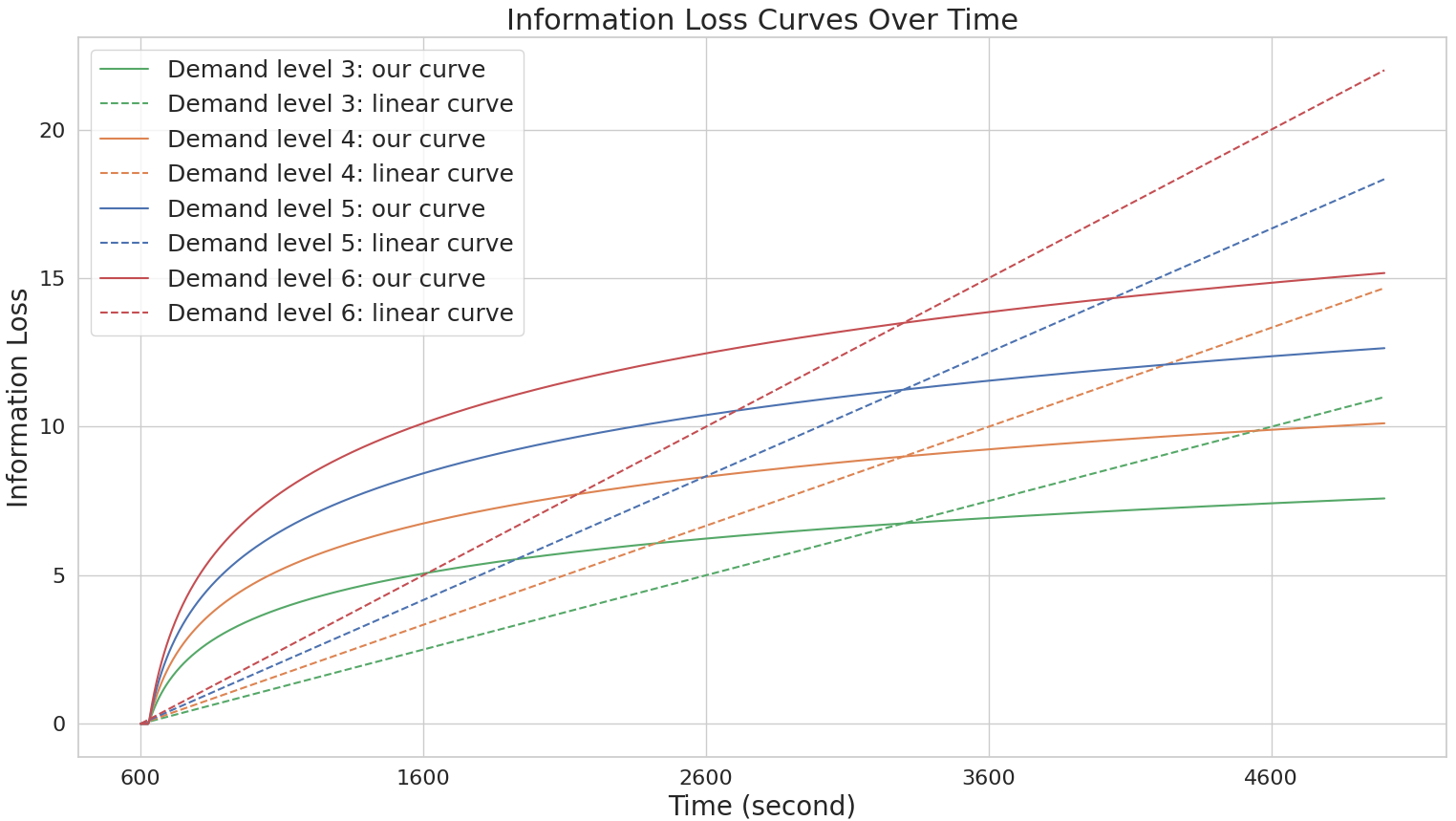}
  \end{center}
  \caption{Comparison of information loss curves.}
  \label{fig:info_compare}
\end{figure}

The comparison between our calibrated information loss curve and the linear information loss curve under different demand levels is shown in Figure \ref{fig:info_compare}.

\subsubsection{Baseline without Information Loss Curve (\textbf{No Info})}
The second baseline model does not incorporate an information loss curve. Instead, in Step 2 of Algorithm \ref{alg:decomp_algorithm}, edge weights are directly updated based on the demand levels from observed data, without accounting for the revisit intervals.

\subsubsection{\textbf{OTOP-RV} Baseline}
To compare our proposed framework with other formulations, we employ an arc-based variant of the Open Team Orienteering Problem with Repeatable Visits (OTOP-RV) model. 
This model, as proposed in reference~\cite{OTOP-RV}, is also designed for urban monitoring tasks that require continuous observation. It extends the traditional TOP by allowing multiple visits to accumulate weighted scores with a constant exponential decay rate.
While our model is designed for the same task, it enables repeated visits only across multiple flights and updates weights dynamically. The similarities and distinctions between our model and the OTOP-RV model make it a suitable baseline for comparison. For a fair comparison, we revise the original formulation of the OTOP-RV model to arc-based representation. Detailed formulations of the OTOP-RV model are provided in the Appendix.

\subsection{Results and Analysis under Static Demand}
We now present and analyze the results of our dynamic multiple-flight framework under the static demand scenario. Initially, all edges are assigned equal weights and demand levels when optimizing the routing plan for the first flight. Using the information loss curve, we generate the optimal routing plan, collect data from observed edges, estimate traffic states for unobserved edges using the most recently observed data or assume zero if no observations are available, and update the edge weights for the next flight. This process is repeated until the maximum number of iteration is reached.

\subsubsection{Two-iteration Routing Results}
\begin{figure*}[!t]
    \centering
    \subfloat[]{\includegraphics[width=0.45\textwidth]{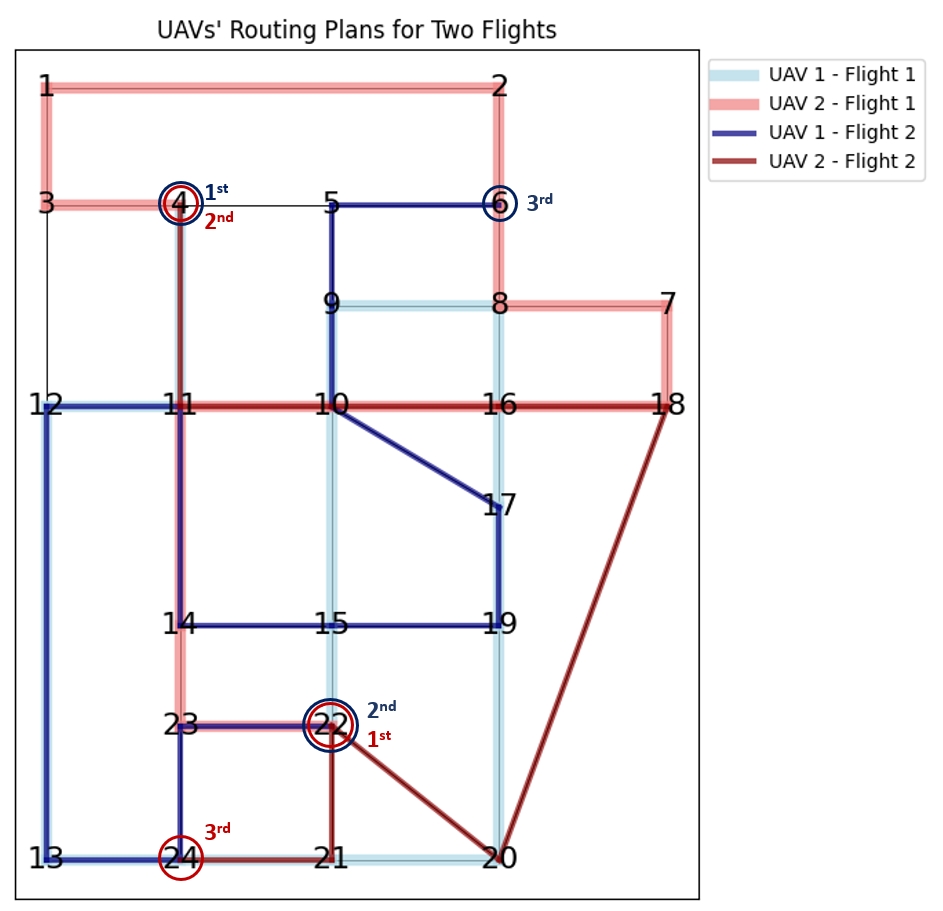}
    \label{fig:single}}
    \hfill
    \subfloat[]{\includegraphics[width=0.45\textwidth]{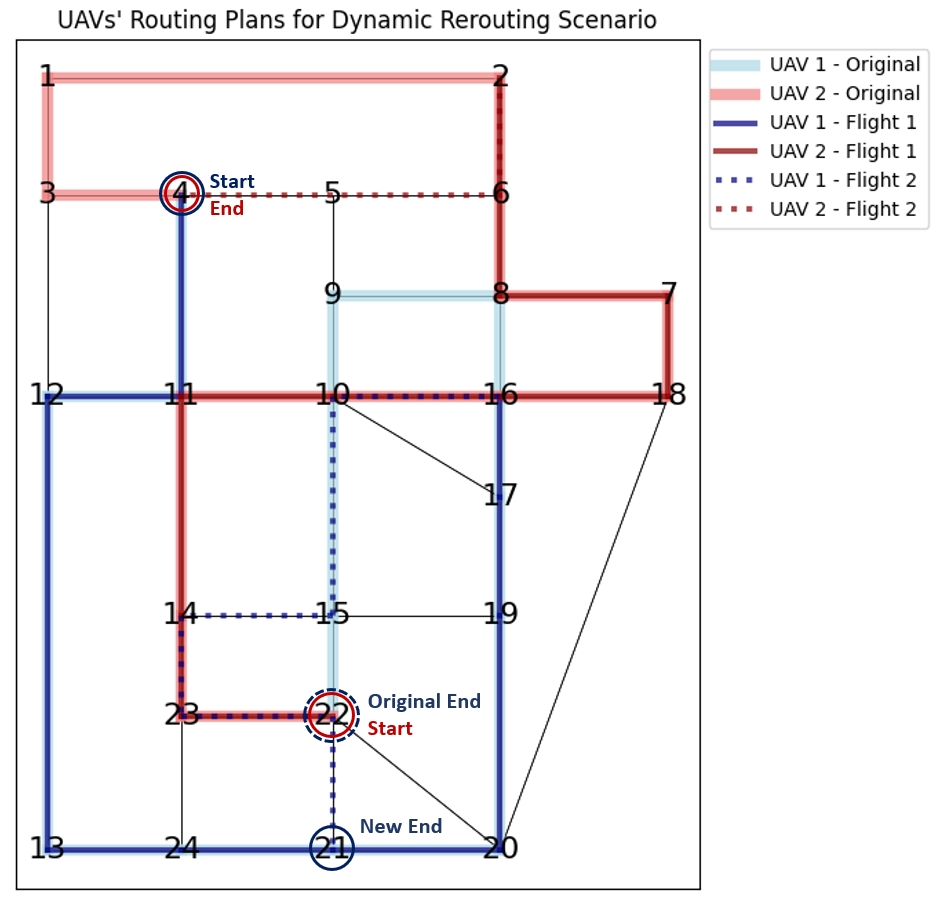}
    \label{fig:single_dynamic}}
    \caption{UAVs' routing results: (a). Routing plans for two flights; (b). Updated routes in a  dynamic rerouting scenario.}
    \label{fig:sioux}
\end{figure*}

To start, we show the results from a simple two-iteration scenario. The flying sequence for UAV 1 is from vertex 4 to vertex 22, then to vertex 6. For UAV 2, the sequence is from vertex 22 to vertex 4, and then to vertex 24. The routing results are illustrated in Figure \ref{fig:single}. As shown in the figure, the two UAVs tend to avoid revisiting edges that have been previously visited, thereby extending the total coverage of the network.

\subsubsection{Dynamic Rerouting Scenario}

Next, we demonstrate a dynamic rerouting scenario where the rendezvous point is changed during a single flight, necessitating the adaptive rerouting. Such changes may arise due to uncertainties in GV operations, such as traffic congestion or unexpected dwell times at bus stops. Specifically, during the flight from vertex 4 to 22 for UAV 1 and from vertex 22 to 4 for UAV 2, the rendezvous point for UAV 1 is changed to vertex 21 instead of vertex 22 after 20 minutes, while UAV 2 continues along its original path to vertex 4. We dynamically reroute the UAVs by treating the completed segment as the first flight, updating edge weights based on the latest traffic data, and re-optimizing the remaining segment as a new flight. The routing results are shown in Figure \ref{fig:single_dynamic}. Although only the rendezvous point for UAV 1 changes, both UAVs adjust their routing plans for the remaining segment to collectively maximize the coverage of the network. 

\subsubsection{Multi-Iteration Performance and Comparative Analysis}

\begin{figure*}[!t]
    \centering
        \subfloat[]{\includegraphics[width=0.45\textwidth]{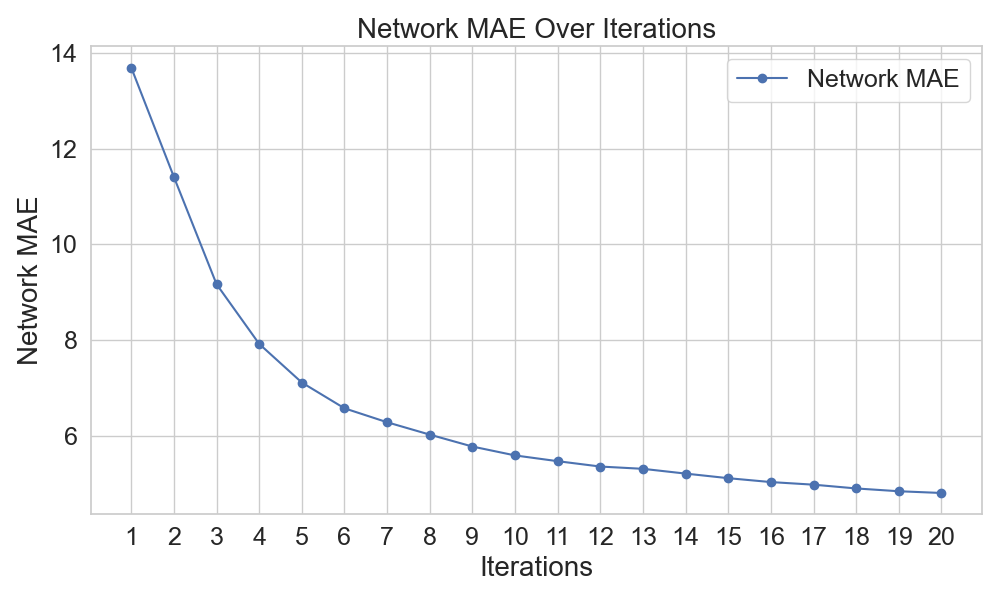}
        \label{fig:iter_mae}}
    \hfill
        \subfloat[]{\includegraphics[width=0.45\textwidth]{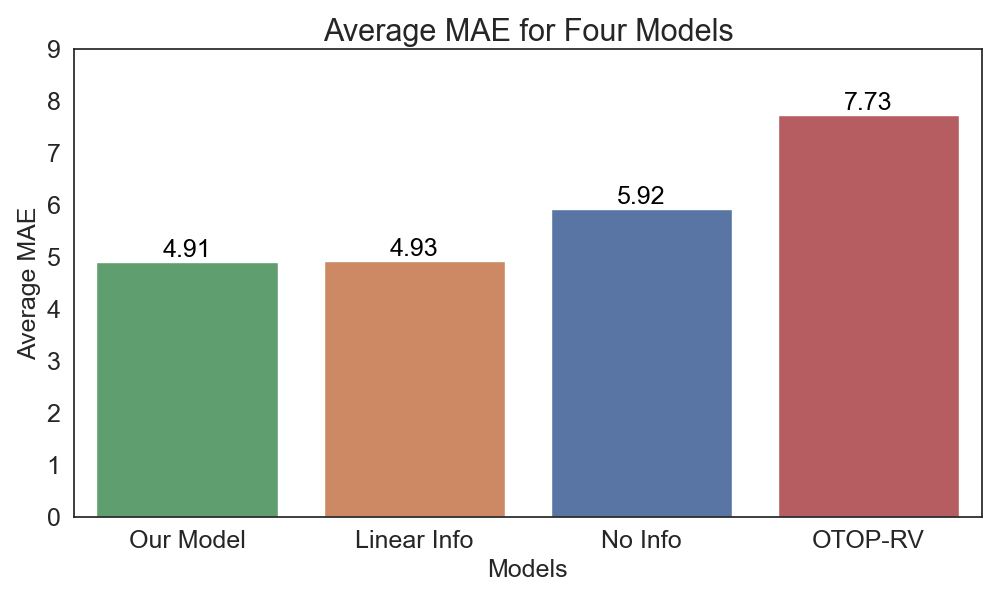}
        \label{fig:bar}}
    \caption{Multi-iteration performance: (a). Network MAE over 20 iterations; (b). Average MAE comparison across four models. }
    \label{fig:multi-result}
\end{figure*}

To evaluate the long-term effectiveness of our dynamic framework, we implement it over 20 iterations. Figure \ref{fig:iter_mae} illustrates the change in network average MAE over all iterations. The estimation error significantly decreases over iterations and eventually stabilizes, demonstrating the framework's ability to improve accuracy over time. As more edges are observed, the weight updates become increasingly precise, and the information loss curve effectively guides UAVs to prioritize edges with higher estimation errors.

To comprehensively evaluate our model's performance, we extract 27 scenarios from the 100-hour dataset, with a 200-minute interval between consecutive scenario start times. Each scenario lasts for 10 hours, allowing us to thoroughly test the 20-iteration case study alongside three baseline models. The overall comparative results are shown in Figure \ref{fig:bar}. Our model achieves an average MAE of 4.91 vehicles, outperforming the OTOP-RV baseline by 36.5\%, the No Info baseline by 17.1\%, and the Linear Info baseline by 0.4\%, demonstrating its superior performance.

\begin{figure*}[!t]
    \centering
        \subfloat[]{\includegraphics[width=\textwidth]{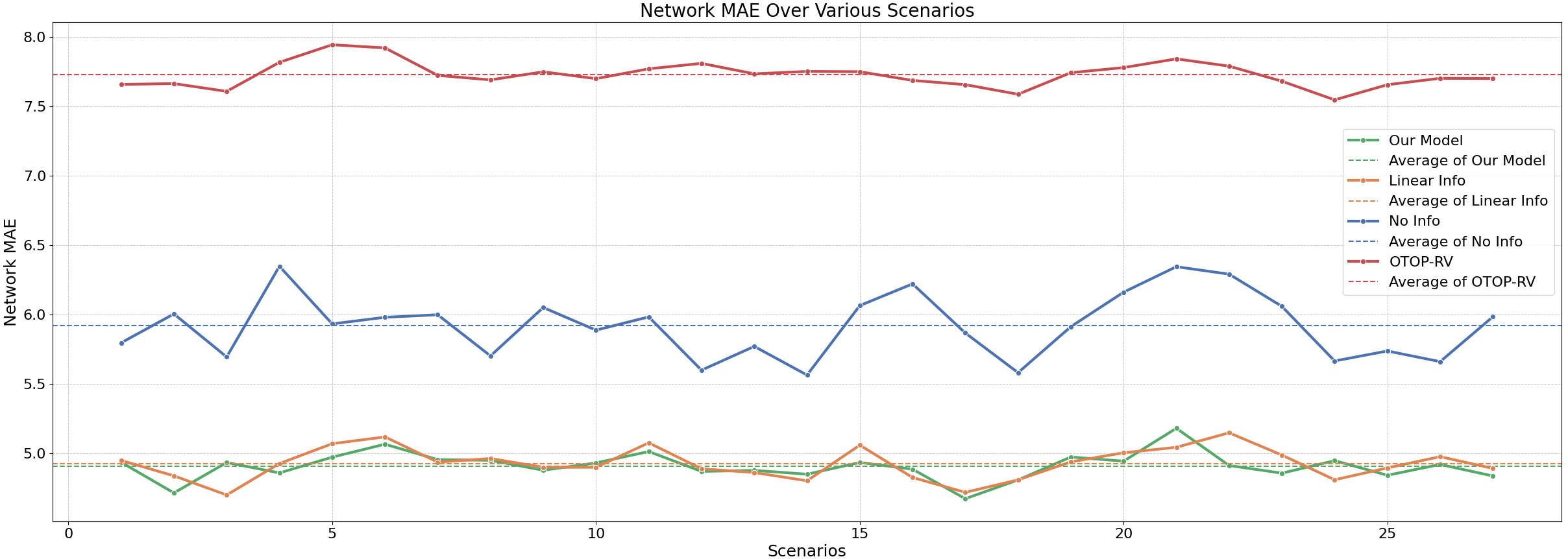}
        \label{fig:TSE0_interval}}
    \hfill
        \subfloat[]{\includegraphics[width=\textwidth]{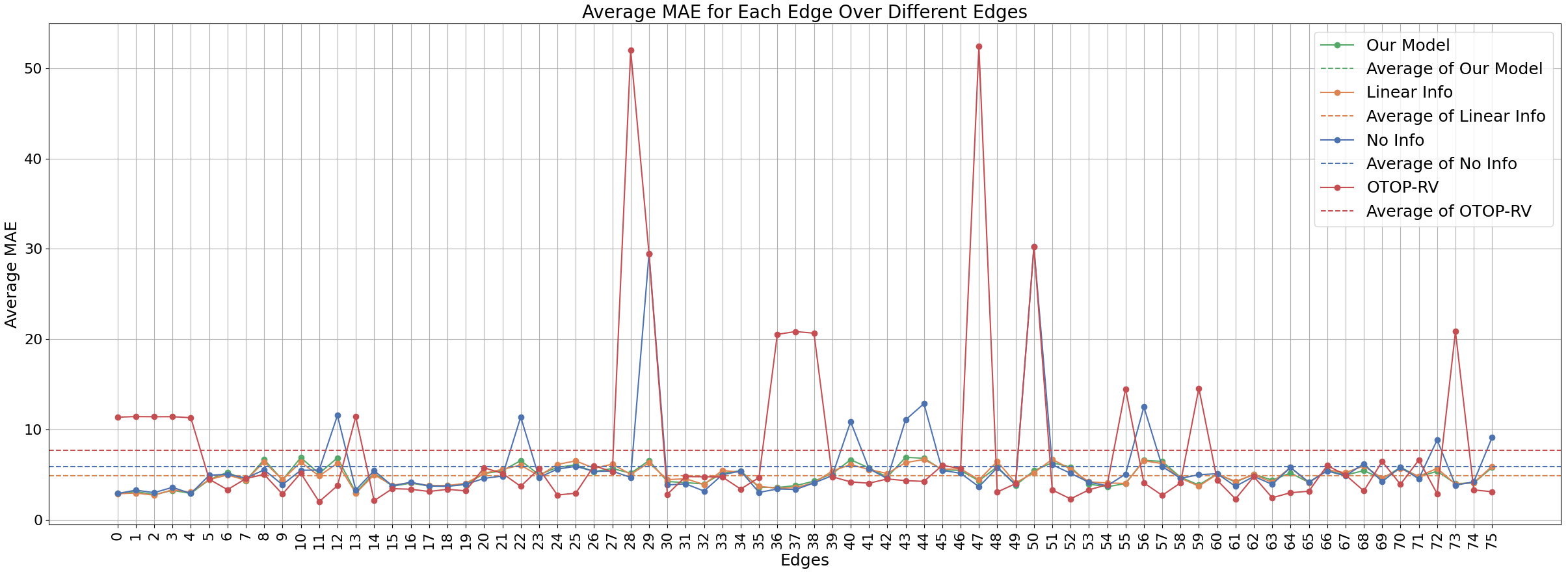}
        \label{fig:TSE0_edge}}
    \caption{Comparative performance for four models: (a). Network MAE over various scenarios; (b). Average MAE for each edge over different edges.}
    \label{fig:TSE0}
\end{figure*}

Figure \ref{fig:TSE0_interval} presents the network MAE across various scenarios of the four models, and Figure \ref{fig:TSE0_edge} illustrates the average MAE of each edge. 
Across all scenarios, the OTOP-RV and No Info baselines consistently show significantly higher errors. In contrast, our model and the Linear Info baseline maintain relatively low errors, consistently outperforming others regardless of the monitoring start time and traffic conditions.
Focusing on edge-level performance, the OTOP-RV baseline shows significant fluctuations, with some edges having errors over 30 vehicles. Similarly, the No Info baseline also shows considerable fluctuations for certain edges. Our model and the Linear Info baseline, however, effectively control the MAE within a narrow range, showing smaller fluctuations and greater stability across the network.

The high errors and fluctuations in the OTOP-RV baseline can be attributed to its static weight and decay functions, which neither integrate real-time traffic conditions nor differentiate between recently visited and long-unvisited edges. This leads to frequent revisits of the same edges, neglecting others for extended periods, and causing large estimation errors due to lack of information. 
Although the No Info baseline considers real-time traffic condition, but it does not account for revisit intervals. It tends to revisit edges with higher flow, ignoring long-unvisited edges with lower flow, resulting in potential inefficiencies and inferior performance compared to our model and the Linear Info baseline.

In contrast, our model and the Linear Info baseline dynamically update weights based on real-time traffic conditions and information loss changes. This enables effective prioritization of high flow edges and long-unvisited edges, maximizing spatiotemporal coverage. The strong performance of both models indicates the effectiveness of the proposed dynamic framework and the design of the weight update mechanism. Moreover, both models demonstrate robustness and stability across the entire network. While our calibrated information loss curve slightly outperforms the linear curve baseline, the overall performance difference is minimal, suggesting that an accurately calibrated information loss curve, though beneficial, may not be essential under our framework, highlighting the scalability of our approach.

\subsection{Sensitivity Analysis on Historical Traffic Information}
In the previous section, we assume no prior traffic information is available for the network. In real-world scenarios, however, historical traffic data (e.g., collected in \cite{pNEUMA}) may be available at certain edges and can be leveraged to estimate traffic demand on unobserved edges. For edges with historical information, historical averages are used to estimate traffic states instead of assuming zero values when prior observations are unavailable. Accordingly, during the optimization of routing plans for the first flight, demand levels are initialized based on the available historical information, while edges lacking such information are set to the lowest demand level. If complete historical information is available, the optimization process begins with the demand levels defined in Equation \ref{eq:demand} and shown in Figure \ref{fig:network}.

\begin{figure}[!ht]
  \begin{center}
    \includegraphics[width=0.8\textwidth]{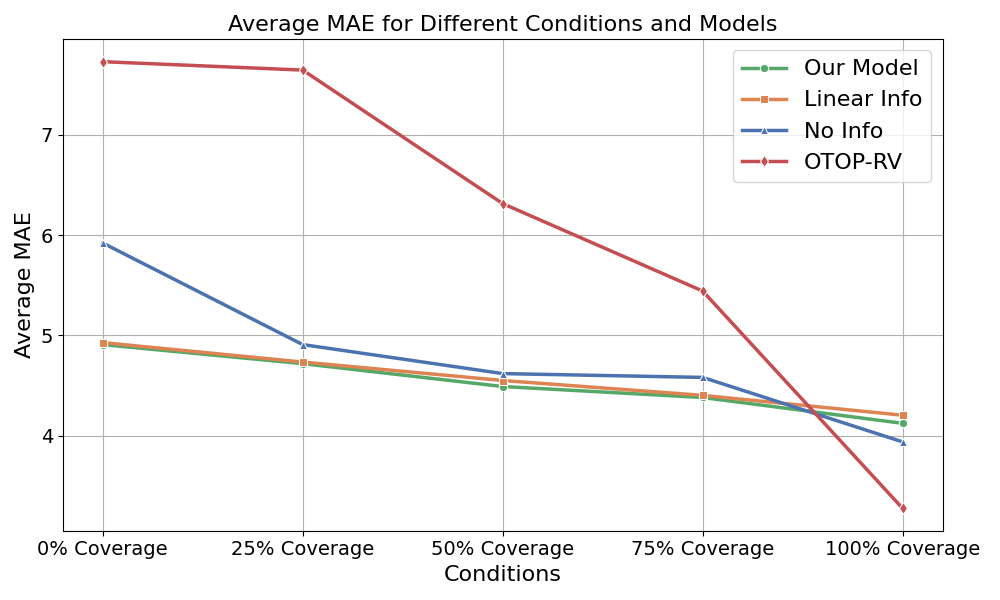}
  \end{center}
  \caption{Average MAE for four models across different coverage conditions.}
  \label{fig:compare}
\end{figure}

We evaluate five conditions: \textbf{0\% Coverage}, \textbf{25\% Coverage}, \textbf{50\% Coverage}, \textbf{75\% Coverage}, and \textbf{100\% Coverage}. For partial coverage, a corresponding percentage of edges is randomly selected. The results are summarized in Figure \ref{fig:compare} and Table \ref{tab:mae}.

\begin{table}[!ht]
    \caption{Average MAE for Four Models under Different Conditions}\label{tab:mae}
    \begin{center}
        \begin{tabular}{l l l l l}
            \textbf{Conditions} & \textbf{Our Model} & \textbf{Linear Info} & \textbf{No Info} & \textbf{OTOP-RV} \\\hline
            \textbf{0\% Coverage} & \textbf{4.908} & 4.928 & 5.921 & 7.729 \\
            \textbf{25\% Coverage} & \textbf{4.718} & 4.734 & 4.908 & 7.645 \\
            \textbf{50\% Coverage} & \textbf{4.490} & 4.550 & 4.619 & 6.314 \\
            \textbf{75\% Coverage} & \textbf{4.381} & 4.402 & 4.581 & 5.442 \\
            \textbf{100\% Coverage} & 4.123 & 4.204 & 3.938 & \textbf{3.275} \\\hline
        \end{tabular}
    \end{center}
\end{table}

Several key insights emerge from the results. First, the estimation accuracy improves across all models as historical information availability increases. Second, both our model and the Linear Info baseline consistently show stable performance across all conditions, with particularly strong performance when information is lacking or incomplete. Third, while the OTOP-RV baseline performs the best with 100\% historical information coverage, it shows significant fluctuations, with average MAE ranging from 3.3 to 7.7 vehicles across all conditions. When historical information is incomplete, even with high coverage (e.g., 75\%), its MAE increases sharply.
This dependency on complete historical data and instability under incomplete information can be attributed to its static framework design, as previously discussed.

In contrast, our model shows remarkable stability, with average MAE ranging from 4.2 to 4.9 vehicles across all conditions. Similarly, the Linear Info baseline and No Info baseline also maintain low variability, indicating less sensitivity to the percentage of available historical information under the proposed dynamic framework.

\subsection{Analysis Under Changing Demand Conditions}
\begin{figure}[!ht]
  \begin{center}
    \includegraphics[width=0.6\textwidth]{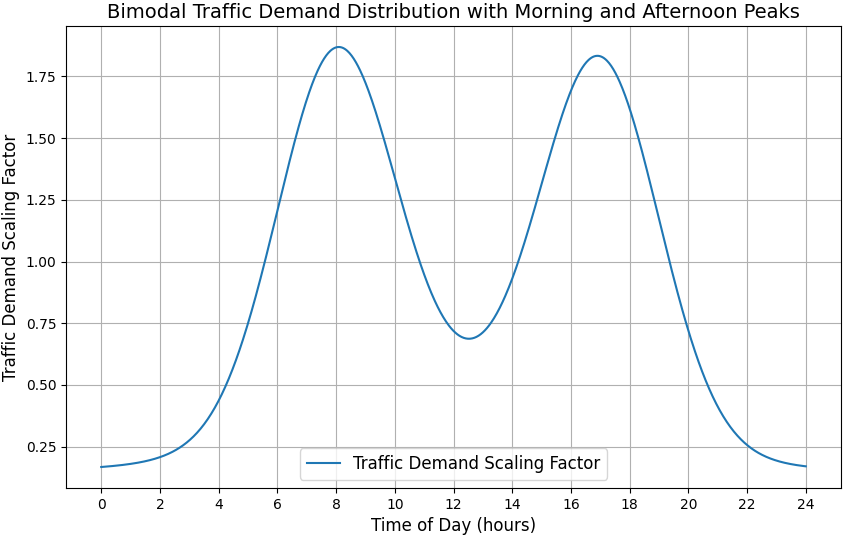}
  \end{center}
  \caption{Changing traffic demand distribution.}
  \label{fig:daily_flow}
\end{figure}

In the previous sections, we analyzed static demand scenarios, where the expected demand remains constant over time. To reflect real-world daily traffic dynamics, we extend our analysis to a changing demand scenario where traffic patterns evolve throughout the day. Inspired by the general shape of daily traffic patterns described in reference~\cite{daily-flow-well-defined-mfd}, we generate a 24-hour traffic demand distribution that captures typical morning and afternoon peaks, as shown in Figure \ref{fig:daily_flow}. A scaling factor is extracted for each hour, which is used to scale the original static demand, creating a 24-hour dataset that mimics real-world traffic evolution. All other simulation setups remain unchanged.

\subsection{Results under Changing Demand}
\begin{figure}[!ht]
  \begin{center}
    \includegraphics[width=0.8\textwidth]{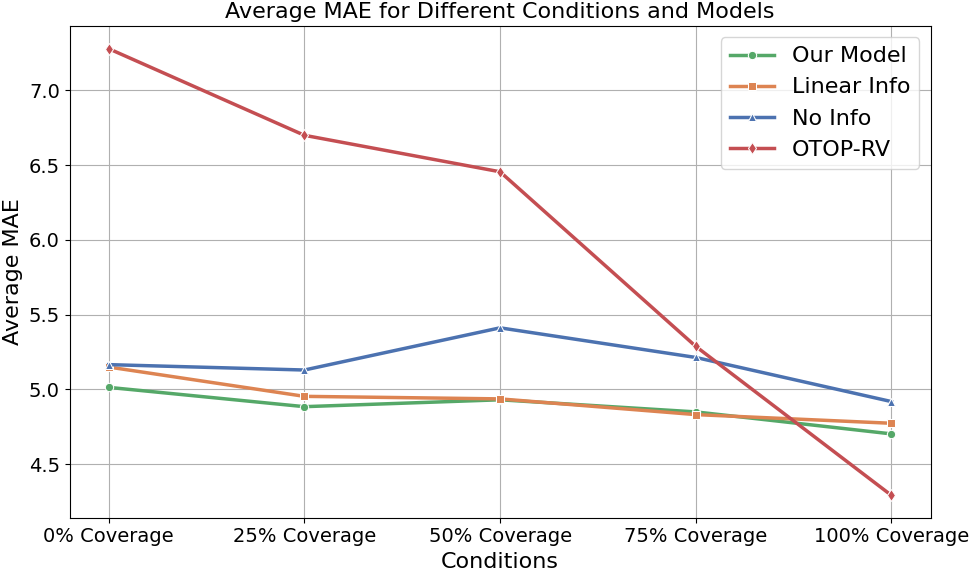}
  \end{center}
  \caption{Average MAE for four models across different coverage conditions under changing demand.}
  \label{fig:compare-changing}
\end{figure}

For the changing demand case study, we analyze 15 scenarios extracted from the 24-hour dataset, with a 50-minute interval between consecutive scenario start times. Each scenario lasts 10 hours, equivalent to 20 iterations. Comparative results under varying historical information coverage are presented in Figure \ref{fig:compare-changing} and Table \ref{tab:mae-changing}.

\begin{table}[!ht]
    \caption{Average MAE for Four Models across Different Conditions under Changing Demand}\label{tab:mae-changing}
    \begin{center}
        \begin{tabular}{l l l l l}
            \textbf{Conditions} & \textbf{Our Model} & \textbf{Linear Info} & \textbf{No Info} & \textbf{OTOP-RV} \\\hline
            \textbf{0\% Coverage} & \textbf{5.014} & 5.150 & 5.166 & 7.278 \\
            \textbf{25\% Coverage} & \textbf{4.885} & 4.954 & 5.130 & 6.699 \\
            \textbf{50\% Coverage} & \textbf{4.931} & 4.937 & 5.412 & 6.455 \\
            \textbf{75\% Coverage} & 4.850 & \textbf{4.832} & 5.214 & 5.288 \\
            \textbf{100\% Coverage} & 4.703 & 4.774 & 4.919 & \textbf{4.292} \\\hline
        \end{tabular}
    \end{center}
\end{table}

The key findings are consistent with those from the static demand scenario. Estimation accuracy improves with increased historical information coverage across all models. The OTOP-RV baseline remains highly sensitive to historical information coverage, exhibiting significant performance variability under lower coverage rates. In contrast, our model, the Linear Info baseline, and the No Info baseline maintain relatively stable performance across all conditions, particularly excelling when information is incomplete or absent, consistently outperforming the OTOP-RV baseline. 
Compared to the static demand scenario, the models utilizing our dynamic framework demonstrate even stronger adaptability to changing demand.

To summarize, our dynamic framework demonstrates adaptability, robustness, and stable performance across all tested conditions including static demand, varying historical information coverage, and changing demand. Our weight update mechanism, combined with the information loss curve design, effectively controls errors within a narrow range across all edges, ensuring consistent monitoring accuracy under varying conditions and throughout the network. This comprehensive approach enables the framework to achieve reliable performance for long-term traffic monitoring tasks under diverse real-world conditions. Its low dependency on historical data and strong adaptivity to dynamic environments make it particularly well-suited for real-world implementations where historical data availability is often limited.

\section{Macroscopic Fundamental Diagram (MFD) Analysis}
\label{sec:mfd}
To achieve a comprehensive and dynamic traffic monitoring system for urban areas, it is essential to understand network-level traffic patterns.
The Macroscopic Fundamental Diagram (MFD) characterizes the relationship between a network's average flow, density, and speed, reflecting the aggregate behavior of traffic within a region~\cite{mfd-1, mfd-2}. 
Despite its potential, few studies have leveraged UAV data to construct the MFD.
While most previous research, such as~\cite{pNEUMA} only constructed link-level fundamental diagrams, a following study from \cite{pNEUMA-2} analyzed multi-modal MFDs during morning peak hours from the pNEUMA dataset. However, these approaches relied on UAVs statically hovering over fixed points with a complete network coverage (i.e., complete information over time and space) and thus introduce high costs and are not scalable. In our dynamic routing framework, only a small portion of links are monitored at the same time, which is equivalent to the number of UAVs (i.e., 2 out of 76 in our case study). The data collected are much more sparse. Therefore, it is interesting to investigate whether the proposed monitoring strategy (i.e., prioritizing high-demand and long time unvisited links) is able to capture the shape and evolution of the MFD with limited data. 

Using the changing demand dataset described in the previous section, we construct the MFDs by computing unweighted network averages of flow, density, and speed at each time interval. Flow is measured as the number of vehicles per 10-minute interval on each edge.
During testing, we assume no prior information is available. For unobserved edges, we estimate the  traffic states using the latest observed data or assume zero if no prior observations exist. 
In addition to the three baseline models introduced earlier, we include a fourth baseline: the OTOP-RV model with 100\% historical information, which outperforms our proposed framework at the link-level. 

\begin{figure*}[!t]
    \centering
        \subfloat[]{\includegraphics[width=0.45\textwidth]{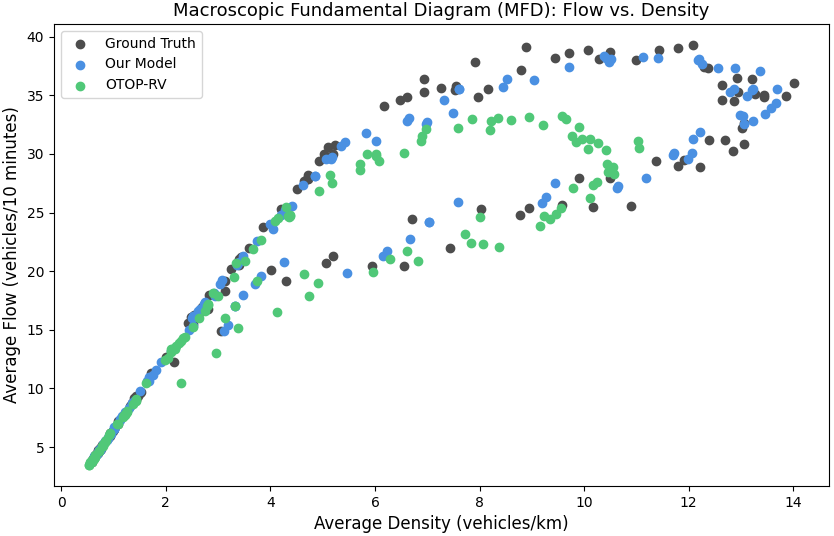} 
        \label{fig:mfd1-sub1}}
    \hfill
        \subfloat[]{\includegraphics[width=0.45\textwidth]{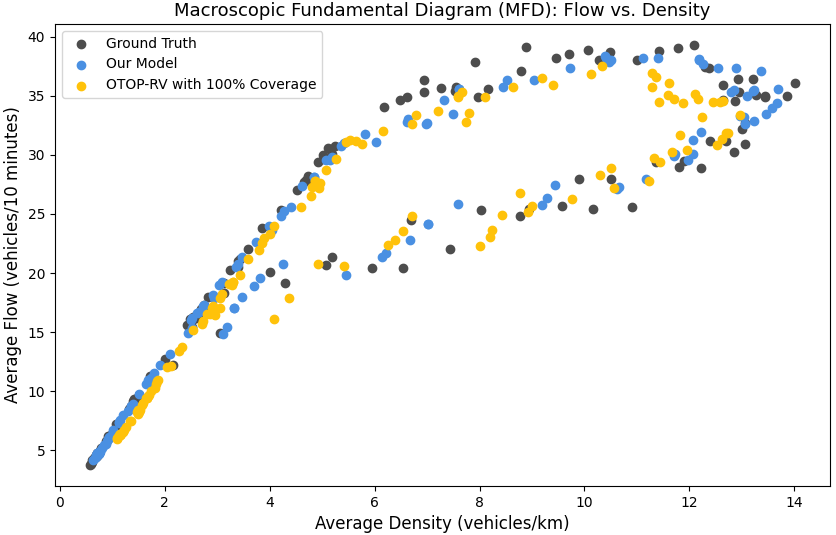} 
        \label{fig:mfd1-sub2}}
    \vskip\baselineskip
        \subfloat[]{\includegraphics[width=0.45\textwidth]{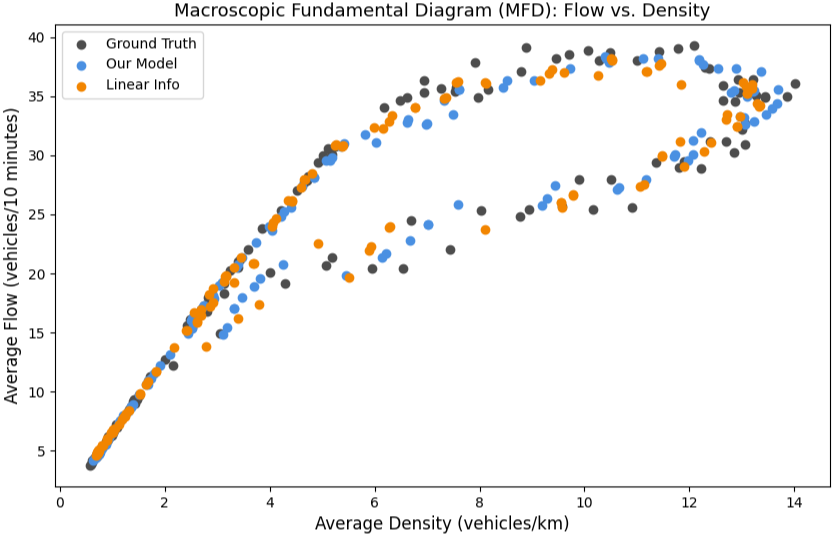} 
        \label{fig:mfd1-sub3}}
    \hfill
        \subfloat[]{\includegraphics[width=0.45\textwidth]{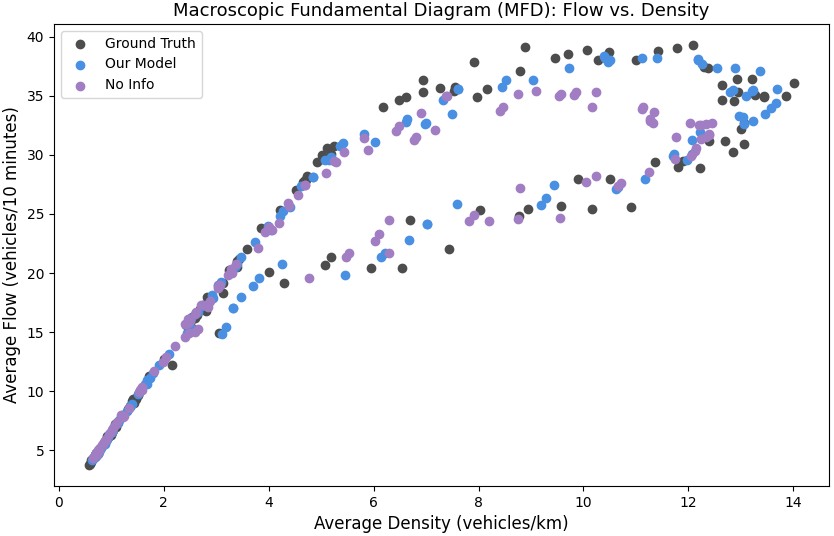} 
        \label{fig:mfd1-sub4}}
    \caption{MFD: Average flow vs. density. Comparison of ground truth, our model's captured MFD, and baseline MFDs.  (a). OTOP-RV baseline; (b). OTOP-RV baseline with 100\% coverage; (c). Linear baseline; (d). No Info baseline.}
    \label{fig:mfd1}
\end{figure*}

\begin{figure*}[!t]
    \centering
        \subfloat[]{\includegraphics[width=0.45\textwidth]{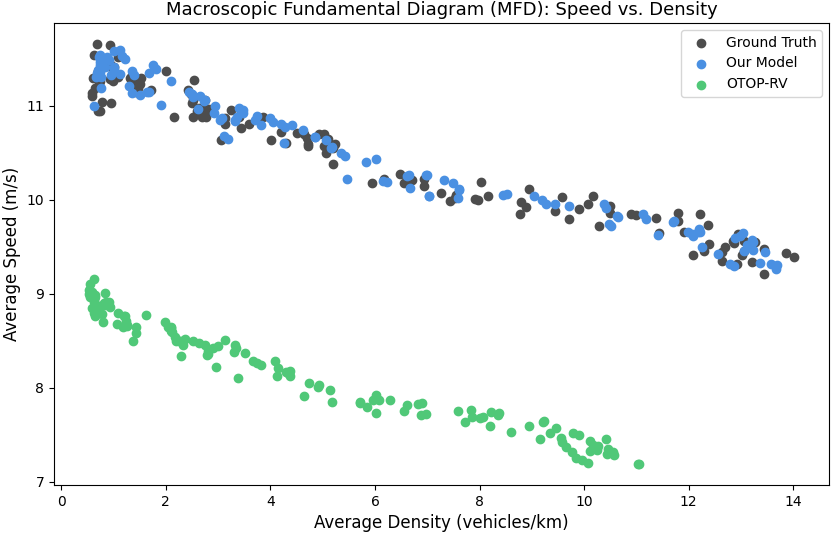} 
        \label{fig:mfd2-sub1}}
    \hfill
        \subfloat[]{\includegraphics[width=0.45\textwidth]{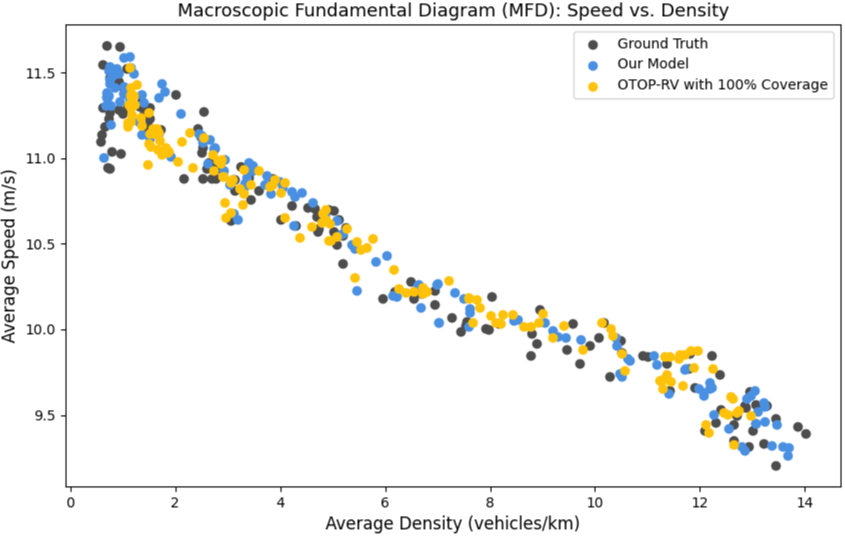} 
        \label{fig:mfd2-sub2}}
    \vskip\baselineskip
        \subfloat[]{\includegraphics[width=0.45\textwidth]{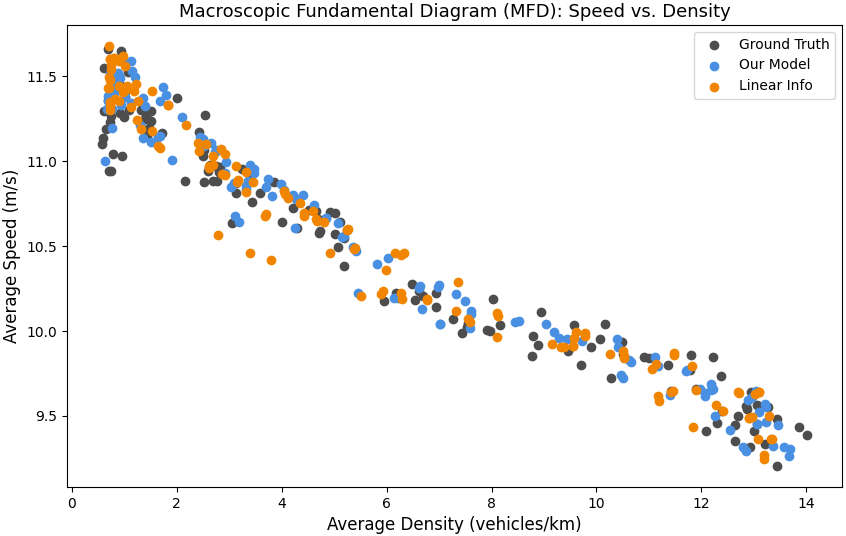} 
        \label{fig:mfd2-sub3}}
    \hfill
        \subfloat[]{\includegraphics[width=0.45\textwidth]{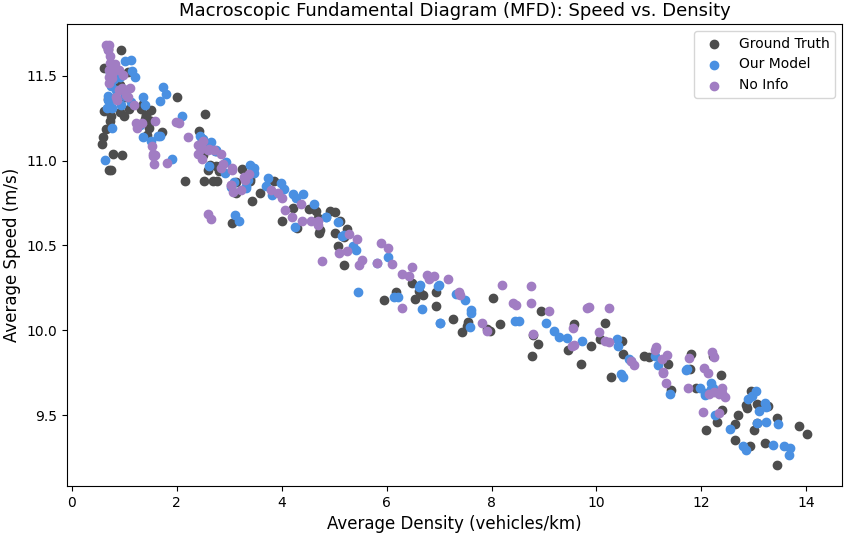}
        \label{fig:mfd2-sub4}}
    \caption{MFD: Average speed vs. density. Comparison of ground truth, our model's captured MFD, and baseline MFDs. (a) OTOP-RV baseline; (b) OTOP-RV baseline with 100\% coverage; (c) Linear baseline; (d) No Info baseline.}
    \label{fig:mfd2}
\end{figure*}

\begin{figure*}[!t]
    \centering
        \subfloat[]{\includegraphics[width=0.45\textwidth]{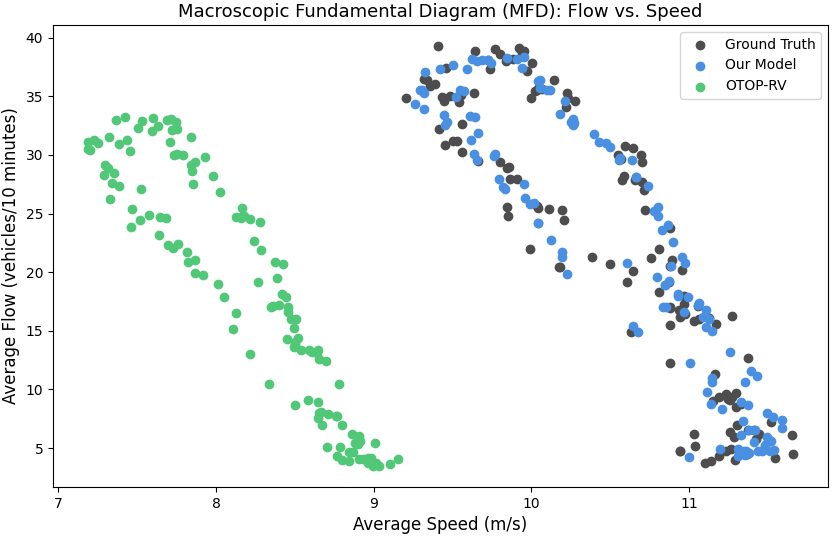} 
        \label{fig:mfd3-sub1}}
    \hfill
        \subfloat[]{\includegraphics[width=0.45\textwidth]{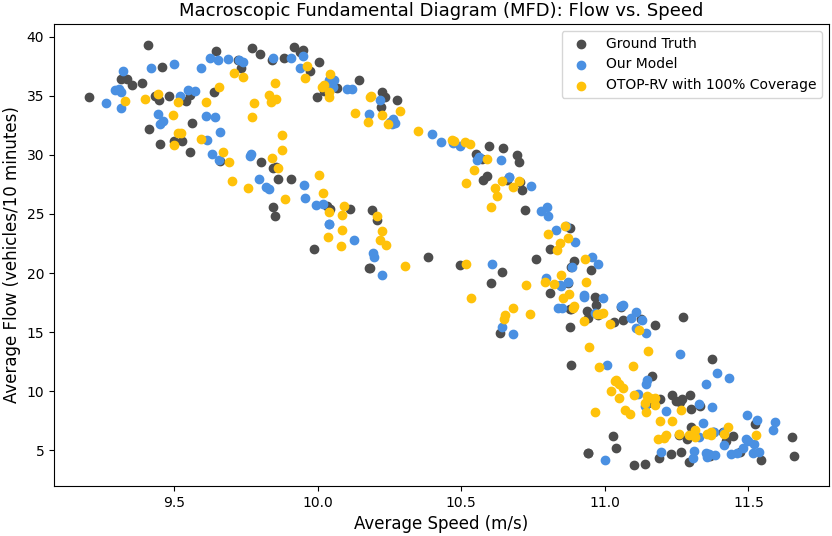} 
        \label{fig:mfd3-sub2}}
    \vskip\baselineskip
        \subfloat[]{\includegraphics[width=0.45\textwidth]{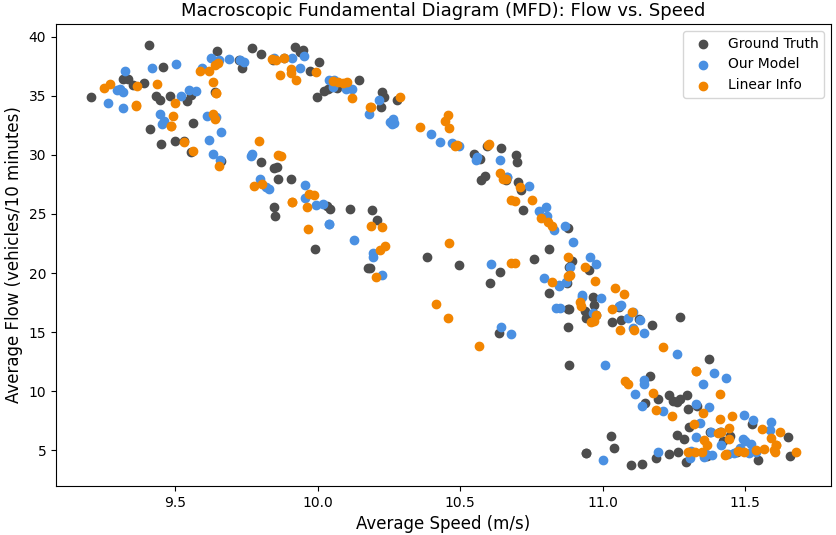} 
        \label{fig:mfd3-sub3}}
    \hfill
        \subfloat[]{\includegraphics[width=0.45\textwidth]{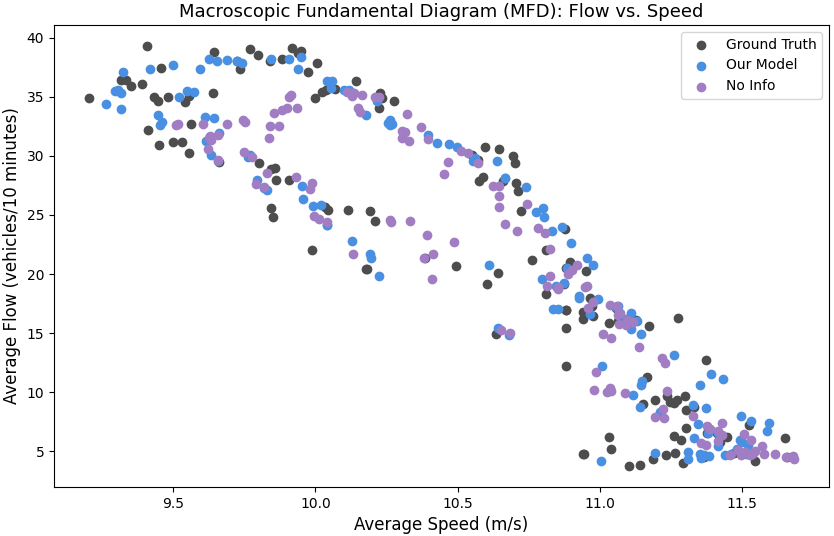}
        \label{fig:mfd3-sub4}}
    \caption{MFD: Average flow vs. speed. Comparison of ground truth, our model's captured MFD, and baseline MFDs. (a) OTOP-RV baseline; (b) OTOP-RV baseline with 100\% coverage; (c) Linear baseline; (d) No Info baseline.}
    \label{fig:mfd3}
\end{figure*}

For each model, we conduct a 24-hour, 48-iteration test to fully capture the long-term network-wide traffic trends as shown in Figure \ref{fig:daily_flow}. We present results for the latter 22 hours, as the demand during the initial two hours is minimal. 
The results are illustrated through three typical MFD plots: average flow versus density (Figure~\ref{fig:mfd1}), speed versus density (Figure~\ref{fig:mfd2}), and average flow versus speed (Figure~\ref{fig:mfd3}). 
Each plot includes the ground truth MFD, the MFD captured by our framework, and those captured by one of the baseline models, allowing for a comparative performance analysis.

The results demonstrate that our dynamic routing framework accurately captures the MFD curves and trends across all three plots.
In the first MFD plot (Figure~\ref{fig:mfd1}), we observe that the OTOP-RV baseline exhibits the greatest deviation from the ground truth curve, while the No Info baseline also shows significant errors in capturing the network-level trend. Even with complete historical information, the OTOP-RV baseline with 100\% coverage (Figure~\ref{fig:mfd1-sub2}) fails to match the accuracy of our framework, highlighting its superior performance in representing network-wide traffic patterns. 
The Linear Info baseline performs comparably to our framework, demonstrating similar accuracy in capturing the MFD at the link-level. 

In the second MFD plot (Figure~\ref{fig:mfd2}), depicting average speed versus density, the trend is nearly linear. The curve generated from the OTOP-RV baseline deviates from the ground truth by approximately 2 m/s at equivalent average densities, which is relatively significant. Our framework and the other baselines can accurately capture this linear trend. The third MFD plot (Figure~\ref{fig:mfd3}) presents similar comparative results, strengthening the observations from the first two plots.

Overall, our dynamic routing framework effectively captures the network-level traffic trends, illustrating its suitability and efficacy for real-world long-term network-wide traffic monitoring tasks. Notably, even without any prior information, our framework can accurately reconstruct the MFD solely using UAV-collected data, outperforming the OTOP-RV baseline model with complete historical information. These findings strongly suggest that our framework provides an effective solution for long-term traffic monitoring under real-world traffic demand patterns.

Moreover, in MFD curve for the flow-density relationship, we observe that the MFD curve is not well-defined, forming a loop-like pattern known as the hysteresis effect. This phenomenon often arises from inhomogeneous demand distributions, traffic instabilities, and variations in driving behavior~\cite{daily-flow-well-defined-mfd, hysteresis-2}. In our case, the hysteresis effect is likely resulted from uneven congestion distribution, as illustrated in Figure \ref{fig:network}, where the central area experiences significantly higher congestion. This localized instability contributes to the network-wide hysteresis effect.

\begin{figure*}[!t]
    \centering
    \subfloat[]{\includegraphics[width=0.45\textwidth]{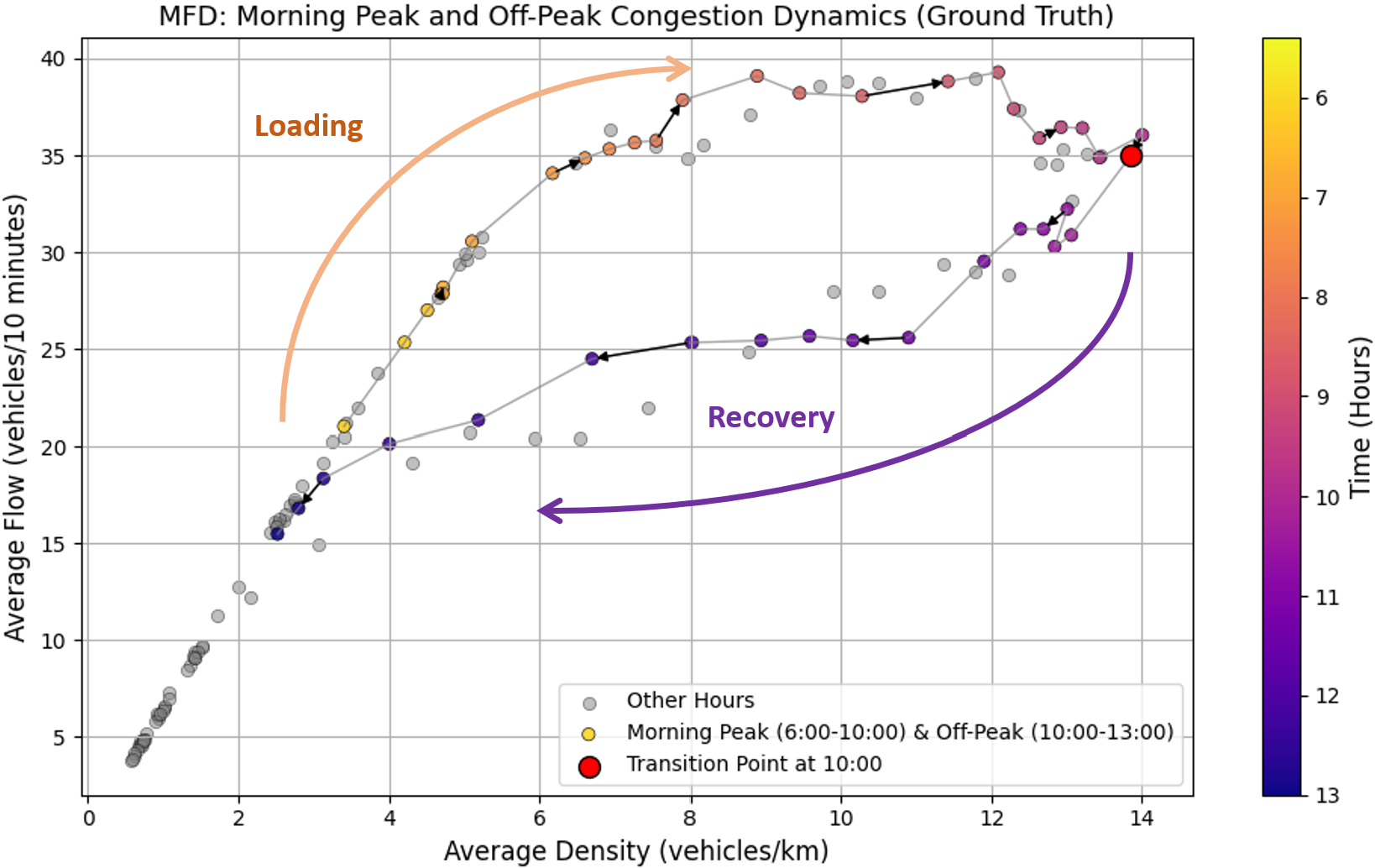} 
        \label{fig:morning-truth}}
    \hfill
    \subfloat[]{\includegraphics[width=0.45\textwidth]{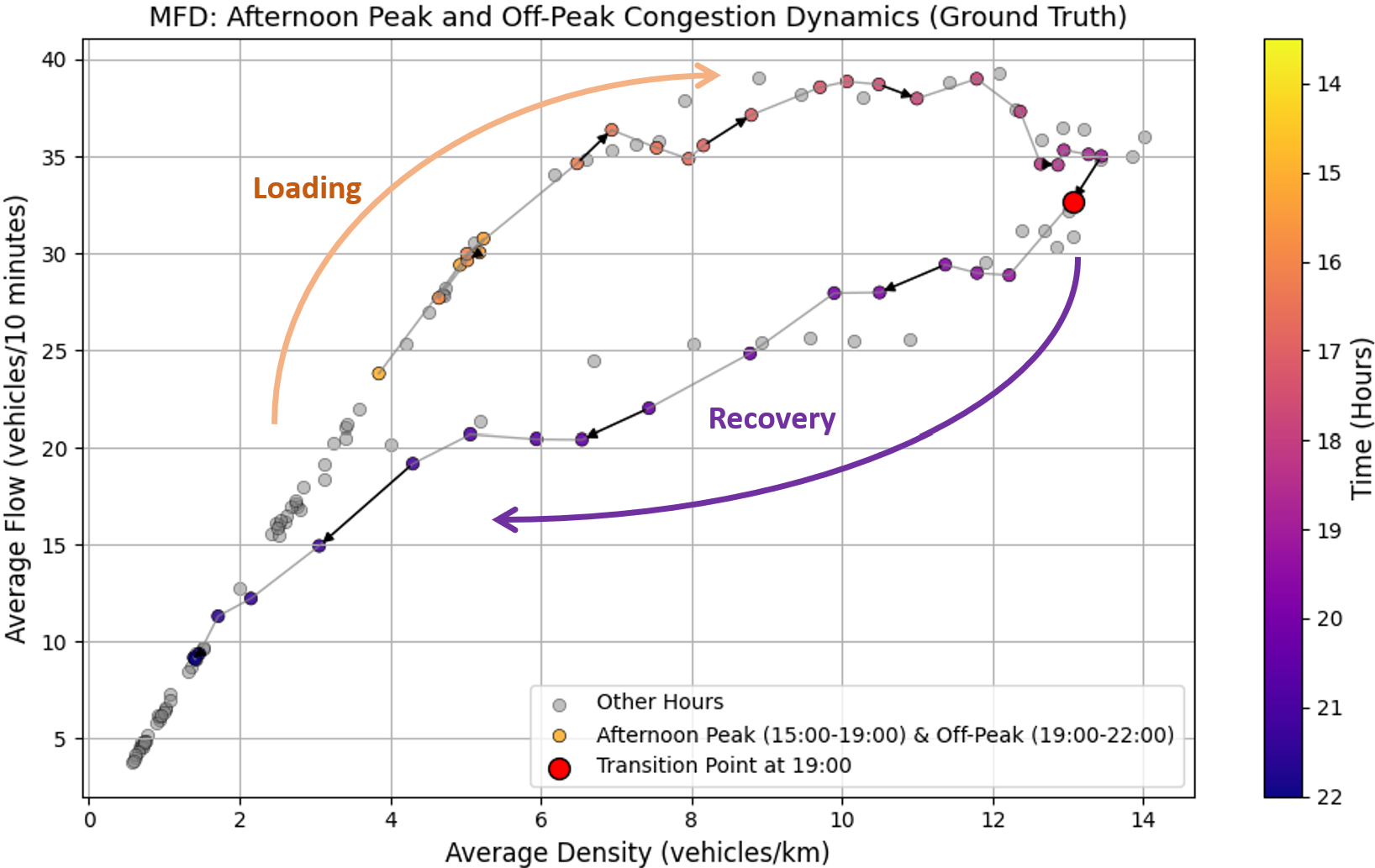}
        \label{fig:afternoon-truth}}
    \vskip\baselineskip
    \subfloat[]{\includegraphics[width=0.45\textwidth]{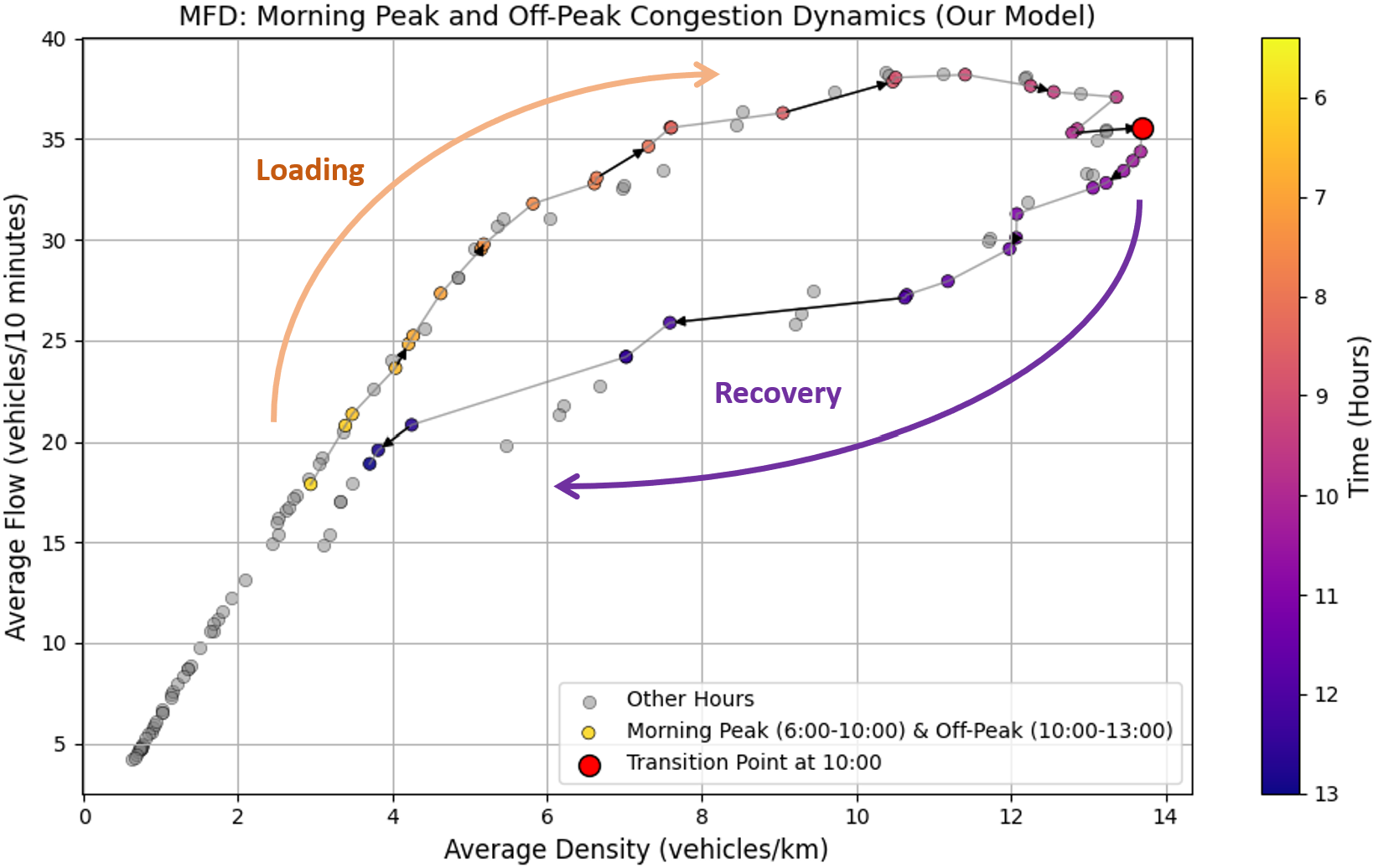}
        \label{fig:morning-our}}
    \hfill
    \subfloat[]{\includegraphics[width=0.45\textwidth]{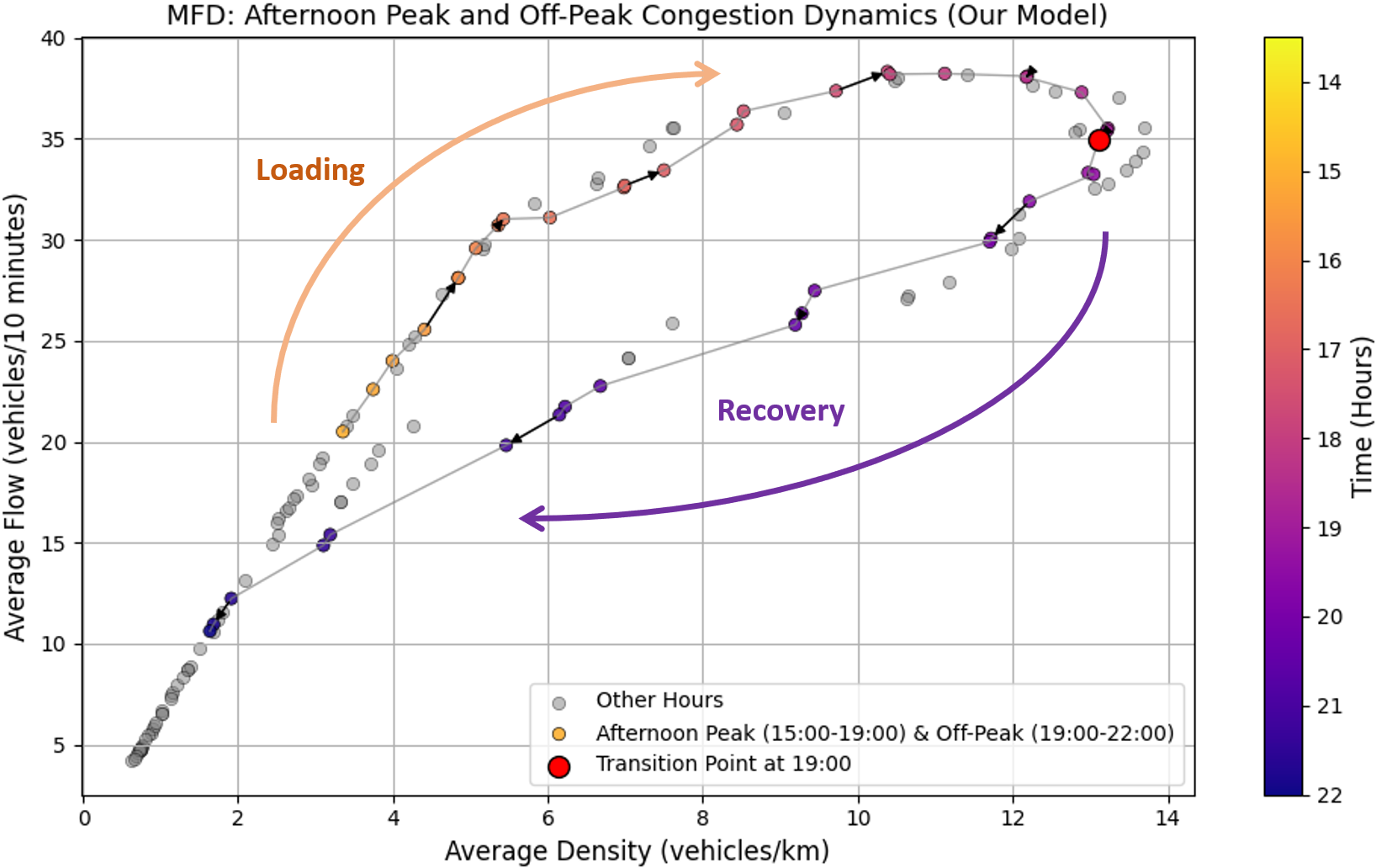}
        \label{fig:afternoon-our}}

    \caption{Comparison of congestion dynamics in MFD hysteresis effects during peak and off-peak hours: (a) Morning peak and off-peak hours in ground truth MFD; (b) Afternoon peak and off-peak hours in ground truth MFD; (c) Morning peak and off-peak hours in constructed MFD; (d) Afternoon peak and off-peak hours in constructed MFD.}
    \label{fig:hysteresis}
\end{figure*}

The hysteresis pattern reflects the temporal evolution of congestion, indicating that congestion does not propagate and dissipate at the same rate~\cite{clockwise-hysteresis}. Figure~\ref{fig:hysteresis} shows the hysteresis effect captured by the ground truth MFD and the MFD constructed by our dynamic routing framework. The color of the data point represents the time it is collected. In the ground truth MFD (Figure~\ref{fig:morning-truth} and Figure~\ref{fig:afternoon-truth}), the clockwise loops align with expectations, showing distinct congestion buildup processes during the morning and afternoon peaks and congestion dissipation processes during mid-day and at night. Specifically, the congestion loading process occurs between 6:00–10:00 for the morning peak and 15:00–19:00 for the afternoon peak, followed by the recovery phases from 10:00–13:00 and 19:00–22:00, respectively. These dynamics, including the two transition points at 10:00 and 19:00, closely correspond to the changing demand distribution used in our simulation settings, as shown in Figure~\ref{fig:daily_flow}.
Similarly, in the MFD captured by our monitoring framework (Figure~\ref{fig:morning-our} and Figure~\ref{fig:afternoon-our}), the temporal trend of congestion evolution is clearly represented, and aligns very well with the ground truth MFD. 
Although there are slight errors around the two transition points (10:00 and 19:00), our framework successfully captures the overall congestion loading and recovery trends, as well as the critical turning points in real time. This ability to capture hysteresis loops in real time is particularly valuable for network-level traffic management, as it reveals evolving congestion dynamics that are typically difficult to access.

\section{DISCUSSION: APPLICATIONS AND FUTURE EXTENSIONS}
\label{sec:discussion}
The main aim of this work is to propose a dynamic and scalable solution for consistent network-level traffic monitoring.
This framework not only intends to collect network-level traffic data more efficiently but also can enable real-time control and management strategies.
In this section, we discuss how the proposed framework can be extended to address GV operation uncertainties, integration with existing sensors, collaboration between UAV and smart GVs (e.g., CAVs), and network-level traffic management strategies informed by the MFD.
We will also discuss potential improvements to the current framework.

\subsection{Incorporating GV Operation Uncertainties}

Our framework leverages predefined GV routes as mobile charging stations for UAVs. In practice, GV travel times can vary due to congestion, incidents, or changes in dwell times for public transit vehicles. To handle such variability, the framework can further incorporate travel time prediction models that use real-time traffic data collected by UAVs. These models can predict future GV locations at the end of UAV flight duration, enabling the adaptive update of rendezvous points. When rendezvous points change, our framework triggers rerouting, terminating the current UAV flight and initiating a new one. This ensures reliable recharging despite uncertainties in GV operations.

\subsection{Real-World Implementation and Data Transfer}

Recent advancements in UAV-based traffic monitoring highlight the feasibility of implementing our framework. UAVs have been successfully used to extract real-time traffic parameters, such as speed, density, and volume, under both free-flow and congested conditions~\cite{real-time-parameter, real-time-parameter-2}, and traffic flow analysis frameworks have been proposed for data collected from moving drones~\cite{moving-UAV}. Leveraging these methods, UAVs in our framework can collect and process traffic data locally during flights and receive route instructions from a pilot or base station as needed. Instead of relying on a fully connected environment, UAVs could offload their monitoring data and receive updated routing information for the next flight during battery swaps at rendezvous points. This approach is cost-effective and technically feasible with current hardware and computer vision advancements. Moreover, integrating additional sensors, such as emission or heat sensors, could further extend the framework's applicability beyond traffic monitoring.

\subsection{Integration with Existing Sensors}
If an urban network already employed fixed or mobile sensors to collect traffic data, our framework can be seamlessly integrated. For instance, when certain links are covered by fixed sensors (e.g., infrastructure-mounted radar/Lidar), their weights can be set to zero, guiding UAVs to focus on other high-priority links. Similarly, if other mobile sensors have recently visited certain edges, we can mark those edges as visited and apply the same weight update mechanism without requiring additional UAV visits. This ensures that the UAV routing strategy remains efficient and avoids duplicating visits where data is already available.

\subsection{CV/CAV-UAV Cooperation}
Beyond traffic monitoring, UAVs can share real-time, network-level information with CV/CAVs, supporting their decision-making process for route selection, lane changing eco-driving and eco-routing, among other operations. This collaboration is particularly valuable because real-time, network-level data is often challenging and costly to acquire, yet it is critical for optimizing CV/CAV performances. Conversely, CV/CAVs equipped with onboard sensors can contribute to the monitoring process while serving as additional mobile charging stations for UAVs. 
In the long term, cooperative CV/CAV-UAV routing could lead to integrated solutions that improve both network-level and link-level traffic efficiency.

\subsection{Real-Time Network-Level Traffic Management with the Generated MFD}
Both link-level and network-level data collected from our framework can be used for MFD-based traffic management strategies. For example, effective perimeter control and route guidance often require partitioning the network into regions with relatively uniform traffic distributions to ensure stable and well-defined MFD curves~\cite{partition,partition-control}. However, creating and maintaining these partitions can be challenging in practice, because real-world traffic patterns are dynamic, spatially diverse, and difficult to capture with complete accuracy.

Our UAV-based approach addresses these challenges by providing accurate and continuous link-level data across the entire network, facilitating the initial partitioning into homogeneous regions suited for MFD estimation. More importantly, because our framework collects data in real time, it enables dynamically adjustments to these partitions as congestion patterns evolve. Instead of relying on fixed regions, real-time network partitioning creates adaptive regional boundaries that reflect current network states.
As a result, our framework can be used to construct real-time regional MFDs and provide continuous updates to perimeter control and route guidance, ensuring management decisions remain responsive and effective.

\subsection{Scalability and Generalization}
Although demonstrated on a small urban network, the proposed framework is inherently scalable and can be applied to larger, more complex networks with challenging traffic conditions. Its adaptability to dynamic, spatiotemporal changes makes it particularly suitable for such scenarios. Beyond urban traffic monitoring, the core concepts of our framework can be generalized to other domains with similar requirements. Potential applications include daily police patrols, emission monitoring, weather surveillance, and disaster response, all of which demand timely, large-scale, and flexible monitoring solutions.

\subsection{Information Loss Curve Improvement}
While our calibrated information loss curve $I(t)$ slightly outperforms the Linear Info baseline with a linear curve, the observed performance difference remains minimal. This suggests that under the current framework, the exact shape of the information loss curve may not critically affect overall monitoring effectiveness. 
Future work should explore whether this limited impact reflects a fundamental property of the framework or if further refinements to the curve could enhance the performance. Insights from information theory may guide the development of more sophisticated curves, potentially improving the framework’s ability to model temporal information loss and capture spatiotemporal patterns more effectively.

\section{CONCLUSIONS}
\label{sec:conclusion}
In this study, we proposed a dynamic UAV routing framework as an efficient and adaptive solution for long-term urban traffic monitoring. We formulated each multi-UAV single-flight routing problem as a TAOP-DP problem, utilizing existing GVs as dynamic charging stations for UAVs without affecting their operations. 
We further decomposed the complex multi-flight problem into manageable single-flight problems, incorporating real-time traffic demand levels and revisit intervals as dynamic weights into the model. 

The framework was tested under various conditions through microscopic simulation on the modified Sioux Falls network. Comparative analyses against three baseline models demonstrated its adaptability, robustness, and stable performance across all experiment settings including static demand, varying historical information, and changing demand, especially in scenarios with incomplete or no prior information.
We also examined the framework's ability to capture the network-wide MFD. Our model could accurately generate the MFD and outperformed all baseline models even when they had access to complete historical information. Its ability to accurately reconstruct network-level traffic patterns like the MFD demonstrates its potential to contribute significantly to network level traffic monitoring and management applications.

\section{APPENDIX}
The Arc-based OTOP-RV formulation is provided below.

\begin{align}
 \max \quad & \sum_{(i,j)\in E} w_{i,j} q_{i,j}^\beta \\
\text{s.t.} \quad & \nonumber\\
& \sum_{i \in V} x_{s_k,i,k} = 1 &\quad \forall k \in K \\
& \sum_{i \in V} x_{i,e_k,k} = 1 &\quad \forall k \in K\\
& \sum_{(i,l) \in E} x_{i,l,k} = \sum_{(l,j) \in E} x_{l,j,k} &\quad \forall l \in V \backslash \{s_k, e_k\}, \forall k \in K\\
& q_{i,j} = 0.5\sum_{k \in K} (x_{i, j, k} + x_{j, i, k})&\quad \forall (i, j) \in E \\
& a_{i,k} + t_{i,j} - a_{j,k} \leq M (1 - x_{i,j,k}) &\quad \forall (i, j) \in E, \forall k \in K\\
& a_{n+1} \leq \Delta \\
& a_{i, k} \geq 0 &\quad \forall i \in V, \forall k \in K\\
& x_{i, j, k} \in\{0,1\}  &\quad \forall i, j \in V, \forall k \in K
\end{align}

The objective function (16) maximizes the weighted rewards within time or distance budget $\Delta$. Constraints (17)-(18) guarantee that each agent departs from its start point $s_k$ and returns to its end point $e_k$. Constraint (19) ensures the flow conservation for each vertex. Constraint (20) counts the number of visitations for each edge, including both directions, and is utilized in the objective function. Constraint (21) prevents the generation of subtours. Constraint (22) ensures that all agents return to their end points within the time or distance budget $\Delta$. Constraints (23)-(24) define the nature of variables $a_{i,k}$ (non negative) and $x_{i,j,k}$ (binary). In this study, we select $\beta = 0.8$.

\bibliographystyle{IEEEtran}
\bibliography{reference}

\vfill

\end{document}